
\magnification=1200
\pageno=1
\baselineskip=18pt

\parindent=0pt

\font\bigbf=cmbx10 scaled\magstep1
\raggedbottom

\input eplain

\newfam\bbbfam
\font\bbbten=msbm10
\font\bbbseven=msbm7
\font\bbbfive=msbm5
\textfont\bbbfam=\bbbten
\scriptfont\bbbfam=\bbbseven
\scriptscriptfont\bbbfam=\bbbfive
\def\bbb{\fam=\bbbfam}

\font\teneufm=eufm10
\font\seveneufm=eufm7
\font\fiveeufm=eufm5
\newfam\eufmfam
\textfont\eufmfam=\teneufm
\scriptfont\eufmfam=\seveneufm
\scriptscriptfont\eufmfam=\fiveeufm

\font\teneufb=eufb10
\font\seveneufb=eufb7
\font\fiveeufb=eufb5
\newfam\eufbfam
\textfont\eufbfam=\teneufb
\scriptfont\eufbfam=\seveneufb
\scriptscriptfont\eufbfam=\fiveeufb

\font\teneurm=eurm10
\font\seveneurm=eurm7
\font\fiveeurm=eurm5
\newfam\eurmfam
\textfont\eurmfam=\teneurm
\scriptfont\eurmfam=\seveneurm
\scriptscriptfont\eurmfam=\fiveeurm

\font\teneurb=eurb10
\font\seveneurb=eurb7
\font\fiveeurb=eurb5
\newfam\eurbfam
\textfont\eurbfam=\teneurb
\scriptfont\eurbfam=\seveneurb
\scriptscriptfont\eurbfam=\fiveeurb

\font\tenscr=rsfs10
\font\sevenscr=rsfs7
\font\fivescr=rsfs5
\skewchar\tenscr='177 \skewchar\sevenscr='177 \skewchar\fivescr='177
\newfam\scrfam \textfont\scrfam=\tenscr \scriptfont\scrfam=\sevenscr
\scriptscriptfont\scrfam=\fivescr
\def\scr#1{{\fam\scrfam#1}}

\rm
\centerline{\bigbf Local square mean in the hyperbolic circle problem }
 \vskip10pt
\centerline{by Andr\'as BIR\'O\footnote{}{Research partially supported by the NKFIH (National Research, Development and Innovation Office) Grants No. K135885, K143876, and by the R\'enyi Int\'ezet Lend\"ulet Automorphic Research Group }}
 \footnote{}{} \footnote{}{2020
Mathematics Subject Classification: 11F72.

Keywords: hyperbolic circle problem, class number of pairs of quadratic forms}
\hfill\break
\centerline{HUN-REN Alfr\'ed R\'enyi Institute of Mathematics}

\centerline{ 1053 Budapest, Re\'altanoda u. 13-15., Hungary; e-mail: biro.andras@renyi.hu}
 \vskip20pt

\noindent {\bf Abstract}. Let $\Gamma\subseteq PSL_2({\bf R})$ be a finite volume
Fuchsian group. The hyperbolic circle problem is the
estimation of the number of elements of the $\Gamma$-orbit of
$z$ in a hyperbolic circle around $w$ of radius $R$, where $z$ and
$w$ are given points of the upper half plane and $R$ is a
large number. An estimate with error term $e^{{2\over 3}R}$ is
known, and this has not been improved for any group.
Petridis and Risager proved that in the special
case $\Gamma =PSL_2({\bf Z})$ taking $z=w$ and averaging over $z$ locally the error
term  can be improved to $e^{\left({7\over {12}}+\epsilon\right)R}$. Here we show such
an improvement for the local $L^2$-norm of the error term.
Our estimate is $e^{\left({9\over {14}}+\epsilon\right)R}$, which is better
than the pointwise bound $e^{{2\over 3}R}$ but weaker than the bound of
Petridis and Risager for the local average.

\vskip20pt

{\bf 1. Introduction.}
\medskip

{\bf 1.1. Statement of the main result.} Let $\bbb H$ be the upper half plane. For $
z,w\in \bbb H$ let
$$u(z,w)={{\left|z-w\right|^2}\over {4\hbox{\rm Im$z$Im$w$}}},\eqno
(1.1)$$
this is closely related to the hyperbolic distance $\rho (z,w)$ of $
z$
and $w$, namely we have $1+2u=\cosh\rho$. The elements of the
group $PSL_2({\bf R})$ act on $\bbb H$ by linear fractional transformations, these are isometries of the
hyperbolic plane. Let $d\mu_z={{dxdy}\over {y^2}}$, this measure is invariant with respect to
the action of $PSL_2({\bf R})$ on $\bbb H$. Let $\Gamma =PSL_2({\bf Z}
)$. For $z\in\bbb H$ and $X>2$ define
$$N\left(z,X\right):=\left|\left\{\gamma\in\Gamma :\hbox{\rm \ $4
u\left(\gamma z,z\right)+2\le X$}\right\}\right|,$$
this is the number of points $\gamma z$ in the hyperbolic circle
around $z$ of radius $\cosh^{-1}\left(X/2\right)$, so the estimation of this
quantity is called the hyperbolic circle problem. We know that
$$\left|N\left(z,X\right)-3X\right|=O_z\left(X^{{2\over 3}}\right
),\eqno (1.2)$$
this is an unpublished theorem of Selberg, but it is
proved also in [L-P], see also [I], Theorem 12.1. Let ${\cal F}$ be
the closure of the standard fundamental domain of $\Gamma$, i.e.
$${\cal F}=\left\{z\in {\bf C}:\;\hbox{\rm Im$z>0$},\>-{1\over 2}
\le\hbox{\rm Re$z$}\le{1\over 2},\>\left|z\right|\ge 1\right\}.\eqno
(1.3)$$
The goal of this paper is to prove the following theorem.

{\bf{\bf THEOREM 1.1. {\it\rm\ {\it Let} $\Gamma =PSL_2({\bf Z})${\it , let} $
{\cal F}$ {\it be as in (1.3), and let} $\Omega\subseteq {\cal F}$ {\it be a compact set, then for any} $
\epsilon >0$ {\it we have}
$$\left(\int_{\Omega}\left(N\left(z,X\right)-3X\right)^2d\mu_z\right
)^{1/2}=O_{\Omega ,\epsilon}\left(X^{{9\over {14}}+\epsilon}\right
).$$
{\bf REMARK 1.1.} The significance of the theorem is that the estimate is better on average than the pointwise bound $
X^{{2\over 3}}$.

{\bf REMARK 1.2.} If $f$ is a smooth nonnegative function which is compactly
supported on ${\cal F}$ and $\epsilon >0$, then
$\int_{{\cal F}}f\left(z\right)\left(N\left(z,X\right)-3X\right)d
\mu_z=O_{f,\epsilon}\left(X^{{7\over {12}}+\epsilon}\right)$ was proved in
[P-R].

{\bf REMARK 1.3.} The bound (1.2) remains valid if we take any finite volume Fuchsian group (a
subgroup of $PSL_2({\bf R})$ acting discontinuously on $\bbb H$ and
having a fundamental domain of finite volume with
respect to $d\mu_z$) in place of $PSL_2({\bf Z})$ provided the main term is
defined including all small Laplace-eigenvalues. Also, the
analogue of the theorem of [P-R] mentioned in Remark 1.2 was proved in [B1] for any finite
volume Fuchsian group with exponent $5/8$ in place of 7/12.

{\bf REMARK 1.4.} It would be interesting to extend Theorem
1.1 for any finite volume Fuchsian group in place of
$PSL_2({\bf Z})$ with some exponent smaller than ${2\over 3}$, similarly as
the theorem of [P-R] was extended in [B1]. Our
present proof uses arithmetic tools, so it might be
extended only for groups similar to $PSL_2({\bf Z})$.

{\bf REMARK 1.5.} Several other kind of average results in the hyperbolic circle problem were proved in [C] and [C-R].

{\bf 1.2. Outline of the proof.} We take an integer $J\ge 2$, it will
be fixed to be large enough in terms of $\epsilon$. We also take a
parameter $d$ which will tend to $\infty$ together with $x$, we
assume $X^{2/3}\le d=X^{1-\delta}$ with some fixed $\delta >0$. We take the sum
$$N_{d,J}\left(z,X\right):=\sum_{j=0}^J\left(-1\right)^j\left(\matrix{
J\cr
j\cr}
\right)\int_1^2\eta_0\left(\tau\right)N\left(z,X-jd\tau\right)d\tau
,$$
where $\eta_0$ is a given nonnegative smooth function on $\left(0
,\infty\right)$ such that $\eta_0\left(\tau\right)=0$
for $\tau\notin\left[1,2\right]$, and $\int_1^2\eta_0\left(\tau\right
)d\tau =1.$ Then the $j=0$ term
equals $N\left(z,X\right)$, but the terms $j\neq 0$ are smoothed versions
of $N\left(z,X\right)$. It can be proved by spectral methods that for
$z\in\Omega$ the $j\neq 0$ terms can be replaced by their main terms
with an error term $O_{\Omega}\left({X\over {\sqrt d}}\right)$. One gets from these
spectral estimates that
$$N_{d,J}\left(z,X\right)=N\left(z,X\right)-3X+O_{\Omega}\left({X\over {\sqrt
d}}\right)\eqno (1.4)$$
for $z\in\Omega$. If we take $d$ larger than $X^{2/3}$, then this error
term will be smaller than $X^{2/3}$.  One can also see easily
that the contribution of the nonhyperbolic $\gamma\in\Gamma$ to $
N_{d,J}\left(z,X\right)$ is
$O_{\Omega ,\epsilon}\left(X^{{1\over 2}+\epsilon}\right)$. Therefore, for the proof of Theorem 1.1 it is enough to estimate
$$\int_{{\cal F}}\left(N_{d,J,{\rm h}{\rm y}{\rm p}}\left(z,X\right
)\right)^2d\mu_z,\eqno (1.5)$$
where $N_{d,J,{\rm h}{\rm y}{\rm p}}\left(z,X\right)$ is the contribution of the hyperbolic
$\gamma\in\Gamma$ to $N_{d,J}\left(z,X\right)$. We will give an expression for (1.5) whose most essential part will be an
expression of type
$$\sum_{t_1,t_2,f^2\neq\left(t_1^2-4\right)\left(t_2^2-4\right)}h\left
(t_1^2-4,t_2^2-4,f\right)\sum_{j_1,j_2=0}^J\left(-1\right)^{j_1+j_
2}\left(\matrix{J\cr
j_1\cr}
\right)\left(\matrix{J\cr
j_2\cr}
\right)F_{X,d}\left(t_1,t_2,f,j_1,j_2\right),\eqno (1.6)$$
where $t_1,t_2>2$ and $f$ run over integers, $F_{X,d}\left(t_1,t_
2,f,j_1,j_2\right)$ is
an analytic expression, and $h\left(t_1^2-4,t_2^2-4,f\right)$ has the
following arithmetic meaning. If $d_1,d_2,t\in {\bf Z}$, then
$h\left(d_1,d_2,t\right)$ denotes the number of $SL_2({\bf Z})$-equivalence
classes of pairs $\left(Q_1,Q_2\right)$ of quadratic forms
$Q_i\left(X,Y\right)=A_iX^2+B_iXY+C_iY^2$ with
integer coefficients satisfying that the discriminant of
$Q_i$ is $d_i$, and the codiscriminant $B_1B_2-2A_1C_2-2A_2C_1$ of
$Q_1$ and $Q_2$ is $t$.

Now, (1.6) can be estimated in the following way. For
certain ranges of the parameters $t_1,t_2$ and $f$ we will
show that if these three parameters are fixed, then the
summation over $j_1,j_2$ will be negligibly small. This will
follow simply from the mean-value theorem of
differential calculus, using that $J$ is large enough. For
those ranges of $t_1,t_2$ and $f$ where this reasoning does
not work, we estimate every term of the summation
separately. In this way we get an upper bound for (1.6) of
size $d^{5/2}X^{-1/2+\epsilon}$. Balancing it with the square of the
error term in (1.4) we
get the theorem choosing $d=X^{5/7}$.

We note that $h\left(d_1,d_2,t\right)$ was studied in the papers [H-W]
and [M]. They gave explicit formulas for $h\left(d_1,d_2,t\right)$ but
only under restrictive conditions for the parameters, so we cannot apply their results. Therefore we prove a general upper bound for $
h\left(d_1,d_2,t\right)$ and apply it in the proof of Theorem 1.1. It would be interesting to investigate in the future whether it is possible to improve the estimate in Theorem 1.1 using an explicit formula instead of our upper bound.

{\bf 1.3. Structure of the paper.} In Section 2 we give a
general formula for the inner product of two automorphic
functions $\sum_{\gamma\in\Gamma_{t_i}}m_i\left(u\left(z,\gamma z\right
)\right)$, where $m_i$ are test
functions, $t_i>2$ are integers for $i=1,2$, and $\Gamma_{t_i}$ is the
set of elements of $SL_2({\bf Z})$ with trace $t_i$. The class
numbers $h\left(t_1^2-4,t_2^2-4,f\right)$ occur in that formula. In
Section 3 we give an upper bound for $h\left(t_1^2-4,t_2^2-4,f\right
)$,
and in Section 4 we investigate the special functions
appearing in the formula of Section 2 in the case when
$m_i$ are characteristic functions as in the circle problem.
In Section 5 we begin the proof of Theorem 1.1 by giving
the spectral estimate and bounding the contribution of
nonhyperbolic elements. In Section 6 we complete the
proof by estimating the square integral (1.5).

{\bf 2. Inner product of automorphic functions and class numbers of pairs of quadratic forms}
\medskip

Our main goal in this section is to prove Lemma 2.2,
which relates the inner product of two automorphic
functions of special kind to class numbers of pairs of
quadratic forms. We prove that lemma in Subsection 2.2.
In Subsection 2.1 we give the necessary definitions and
we prove an easy lemma which will be used later.

{\bf 2.1. The necessary definitions and an upper bound.} We start by taking a positive discriminant $
s$ and
introducing the set $\scr{Q}_s$ of quadratic forms with
discriminant $s$. Let $s$ be a positive integer with $s\equiv 0,1$
modulo $4$ and define
$$\scr{Q}_s:=\left\{Q\left(X,Y\right)=AX^2+BXY+CY^2:\;A,B,C\in {\bf Z},\;
B^2-4AC=s\right\}.$$
If $\tau =\left(\matrix{a&b\cr
c&d\cr}
\right)\in SL_2({\bf Z})$ and $Q$ is a quadratic form, let us
define the quadratic form $Q^{\tau}$ by $Q^{\tau}\left(X,Y\right)
=Q\left(aX+bY,cX+dY\right)$. The
group $SL_2({\bf Z})$ acts in this way on $\scr{Q}_s$. When $s=t^2-4$,
the set $\scr{Q}_s$ can be identified with elements in $SL_2({\bf Z})$
with trace $t$. Indeed, if $t>2$ is an integer, let
$$\Gamma_t=\left\{\left(\matrix{a&b\cr
c&d\cr}
\right)\in SL_2({\bf Z}):\;\;a+d=t\right\}.$$
The group $SL_2({\bf Z})$ acts on this set by conjugation. If
$\gamma =$$\left(\matrix{a&b\cr
c&d\cr}
\right)\in\Gamma_{}$$_t$, let $Q_{\gamma}(X,Y)=cX^2+(d-a)XY-bY^2$.
Then it is easy to see (see [B2], p. 119) that the map
$\gamma\mapsto Q_{\gamma}$ is a one-to-one correspondence between $
\Gamma_t$ and $\scr{Q}_s$ with
$s=t^2-4$, and also between the conjugacy classes of $\Gamma_t$ over $
SL_2({\bf Z})$ and the
$SL_2({\bf Z})$-equivalence classes of $\scr{Q}_s$. More precisely: if
$\tau\in SL_2({\bf Z})$ and $\gamma\in\Gamma_t$$ $, then we have $
Q_{\tau^{-1}\gamma\tau}=Q_{\gamma}^{\tau}$. Observe
also that the fixed points of $\gamma$ on ${\bf R}$ are exactly the
roots of the quadratic polynomial $Q_{\gamma}\left(X,1\right)$.

For $d_1,d_2,t\in {\bf Z}$, let ${\cal Q}_{d_1,d_2,t}$ be the subset of $\scr{Q}_{d_1}\times\scr{Q}_{d_2}$ consisting of those pairs $\left
(Q_1,Q_2\right)$ of quadratic forms having codiscriminant $t$. In other words, writing
$$Q_1\left(X,Y\right)=A_1X^2+B_1XY+C_1Y^2,\;Q_2\left(X,Y\right)=A_
2X^2+B_2XY+C_2Y^2\eqno (2.1)$$
we require that the discriminant of $Q_j$ is $d_j$ ($j=1,2$) and that
$$B_1B_2-2A_1C_2-2A_2C_1=t.\eqno (2.2)$$
It is easy to check that if $\tau\in SL_2({\bf Z})$, and
$\left(Q_1,Q_2\right)\in {\cal Q}_{d_1,d_2,t}$, then $\left(Q_1^{\tau},Q_
2^{\tau}\right)\in {\cal Q}_{d_1,d_2,t}$. Hence $SL_2({\bf Z})$
acts on ${\cal Q}_{d_1,d_2,t}$. Let us denote by $h\left(d_1,d_2,t\right
)$ the number of $SL_2({\bf Z})$-equivalence classes of ${\cal Q}_{d_1,d_
2,t}$.

If $t_1>2$, $t_2>2$ are integers, let ${\cal R}_{t_1,t_2}$ be the subset of $\scr{
Q}_{t_1^2-4}\times\scr{Q}_{t_2^2-4}$ consisting of those pairs $\left
(Q_1,Q_2\right)$ of quadratic forms satisfying that
$$Q_1=\lambda Q_2\hbox{\rm \ with some $\lambda\in {\bf Q}$}.\eqno
(2.3)$$
Note that ${\cal R}_{t_1,t_2}$ is empty unless ${{t_1^2-4}\over {
t_2^2-4}}\in {\bf Q}^2.$ It is easy to check that if $\tau\in SL_
2({\bf Z})$, and
$\left(Q_1,Q_2\right)\in {\cal R}_{t_1,t_2}$, then $\left(Q_1^{\tau},Q_2^{
\tau}\right)\in {\cal R}_{t_1,t_2}$. Hence $SL_2({\bf Z})$
acts on ${\cal R}_{t_1,t_2}$. Let ${\cal R}_{t_1,t_2}^{\ast}$ denote a complete set of representatives of the $SL_2({\bf Z})$-equivalence classes of ${\cal R}_{t_1,t_2}$.

If $\left(Q_1,Q_2\right)\in {\cal R}_{t_1,t_2}$, we can define a
nonnegative real number $n\left(Q_1,Q_2\right)$ in the following way. Using the bijection $\gamma\mapsto Q_{\gamma}$ defined earlier, let
$\gamma_i\in\Gamma_{t_i}$ be such that $Q_{\gamma_i}=Q_i$ for $i=
1,2$. Then $\gamma_i$ are
uniquely determined. The fixed points on ${\bf R}$ of the
hyperbolic transformations $\gamma_1$ and $\gamma_2$ are the same by
(2.3), since they are the roots of the
polynomial $Q_1\left(X,1\right)=\lambda Q_2\left(X,1\right)$. Denoting the centralizer
of a hyperbolic element $\gamma$ in $SL_2({\bf Z})$ by $C\left(\gamma\right
)$ it is well-known and easily proved that we have
$$C\left(\gamma\right)=\left\{\tau\in SL_2({\bf Z}):\;\;\tau z_1=
z_1,\tau z_2=z_2\right\},\eqno (2.4)$$
where $z_1$ and $z_2$ are the fixed points of $\gamma$. Therefore we
have $C\left(\gamma_1\right)=C\left(\gamma_2\right)$. The image of $
C\left(\gamma_1\right)$ in $PSL_2({\bf Z})$ is infinite cyclic, i.e there is a $
\gamma_0\in SL_2({\bf Z})$ such that
$$C\left(\gamma_1\right)=\left\{\pm\gamma_0^l\in SL_2({\bf Z}):\;\;
l\in {\bf Z}\right\}.\eqno (2.5)$$
Let $N\left(\gamma\right)$ denote the norm of a hyperbolic transformation
$\gamma$, see p. 19 of [I]. Let us define $n\left(Q_1,Q_2\right):
=\left|\log N\left(\gamma_0\right)\right|$, this quantity is well-defined. It can be seen that if $\tau\in SL_2({\bf Z})$, then
$n\left(Q_1^{\tau},Q_2^{\tau}\right)=n\left(Q_1,Q_2\right).$ Finally, if $
t_1>2$, $t_2>2$ are integers,
let us define
$$E_{t_1,t_2}:=\sum_{\left(Q_1,Q_2\right)\in {\cal R}_{t_1,t_2}^{
\ast}}n\left(Q_1,Q_2\right).\eqno (2.6)$$
The following lemma will be enough for handling $E_{t_1,t_2}$
during the proof of Theorem 1.1.

{\bf LEMMA 2.1.} {\it If} $2<t_1\le t_2$ {\it are integers, then} $
E_{t_1,t_2}\ll_{\epsilon}t_2^{1+\epsilon}$
{\it for every} $\epsilon >0${\it .}

{\it Proof.\/} If $\left(Q_1,Q_2\right)\in {\cal R}_{t_1,t_2}$ and $
\gamma_i\in\Gamma_{t_i}$ is such that
$Q_{\gamma_i}=Q_i$ for $i=1,2$, then for $\gamma_0$ satisfying (2.5) we
clearly have $\left|\log N\left(\gamma_0\right)\right|\le\left|\log
N\left(\gamma_2\right)\right|\ll\log t_2$. So it is
enough to show that the number of $SL_2({\bf Z})$-equivalence
classes of ${\cal R}_{t_1,t_2}$ is $\ll_{\epsilon}t_2^{1+\epsilon}$. If $
Q_2\in\scr{Q}_{t_2^2-4}$ is given, then
there are at most two possibilities for $Q_1$ to have
$\left(Q_1,Q_2\right)\in {\cal R}_{t_1,t_2}$, so it is enough to show that the number
of $SL_2({\bf Z})$-equivalence classes of $\scr{Q}_{t_2^2-4}$ is $
\ll_{\epsilon}t_2^{1+\epsilon}$. But
this is a well-known statement, it follows at once from
[Bu], Proposition 3.3 and formula (3.1). The lemma is
proved.

{\bf 2.2. The formula for the inner product.} If $t>2$ is an integer and $
m$ is a compactly supported bounded function on
$[0,\infty )$, then for $z,w\in\bbb H$ write
$$m(z,w)=m\left(u(z,w)\right)\eqno (2.7)$$
by an abuse of notation, see (1.1) for $u(z,w)$. For $z\in\bbb H$ define
$$M_{t,m}(z)=\sum_{\gamma\in\Gamma_t}m\left(z,\gamma z\right).\eqno
(2.8)$$
The main result of this subsection, Lemma 2.2 expresses
the inner product of two such functions $M_{t_1,m_1}$, $M_{t_2,m_
2}$
in terms of the quantities $E_{t_1,t_2}$ and $h\left(t_1^2-4,t_2^
2-4,f\right)$
defined above. We need also the following definitions to
state the lemma.

Let $t_1,$$t_2>2$ be real numbers and let $m_1$, $m_2$ be compactly
supported bounded functions on $[0,\infty )$. Let us write
$${\cal J}\left(t_1,t_2,m_1,m_2\right):=\int_{-\pi /2}^{\pi /2}m_
1\left({{t_1^2-4}\over {4\cos^2\theta}}\right)m_2\left({{t_2^2-4}\over {
4\cos^2\theta}}\right){{d\theta}\over {\cos^2\theta}},\eqno (2.9)$$
and for every real $F$ with $\left|F\right|\neq 1$ let us write
$${\cal I}\left(t_1,t_2,F,m_1,m_2\right):=\int\!\!\!\int{{m_1\left
({{t_1^2-4}\over 4}\left(1+S^2\right)\right)m_2\left({{t_2^2-4}\over
4}\left(1+T^2\right)\right)}\over {\sqrt {S^2+T^2+2FTS+1-F^2}}}dS
dT,\eqno (2.10)$$
where we integrate over the set
$$\left\{\left(S,T\right)\in {\bf R}^2:\;S^2+T^2+2FTS+1-F^2>0\right
\}.\eqno (2.11)$$
{\bf REMARK 2.1.} The integral (2.9) is absolutely convergent, because $
m_1$ and $m_2$ are
compactly supported and bounded. The absolute
convergence of the integral in (2.10) is trivial in the case $\left
|F\right|<1$, because then we always have
$S^2+T^2+2FTS\ge 0$. In the case $\left|F\right|>1$ we use the linear
substitution $u=S+\left(F+\sqrt {F^2-1}\right)T$, $v=S+\left(F-\sqrt {
F^2-1}\right)T$.
Since $m_1$ and $m_2$ are compactly supported, we have in
(2.10) that $\left|S\right|$ and $\left|T\right|$ are bounded from above. But then
we have also $\left|u\right|$$,\left|v\right|<C$ with some $C>0$. The condition
in (2.11) reads as $uv\ge F^2-1$, therefore we have also $\left|u\right
|$,$\left|v\right|>c$
with some $c>0$. Since integrating over the set defined
by the conditions $c<\left|u\right|,\left|v\right|<C$, $uv\ge F^2
-1$ we clearly
have $\int\!\!\!\int{1\over {\sqrt {uv+1-F^2}}}dudv<\infty$, we get that (2.10) is absolutely convergent also for $\left
|F\right|>1$.

{\bf LEMMA 2.2.} {\it Let} $t_1,$$t_2>2$ {\it be integers and let} $
m_1${\it ,} $m_2$ {\it be compactly supported bounded functions on} $[0,\infty )${\it . Then using the notation (2.8) we have that }
$$\int_{{\cal F}}M_{t_1,m_1}(z)M_{t_2,m_2}(z)d\mu_z\eqno (2.12)$$
{\it equals the sum of }
$${\cal J}\left(t_1,t_2,m_1,m_2\right)E_{t_1,t_2}\eqno (2.13)$$
{\it and}
$$\sum_{f\in {\bf Z},f^2\neq\left(t_1^2-4\right)\left(t_2^2-4\right
)}h\left(t_1^2-4,t_2^2-4,f\right){\cal I}\left(t_1,t_2,{f\over {\sqrt {
t_1^2-4}\sqrt {t_2^2-4}}},m_1,m_2\right).\eqno (2.14)$$
{\it See (2.1), (2.2), and the paragraph below (2.2) for}
$h\left(t_1^2-4,t_2^2-4,f\right)${\it, (2.6) for} $E_{t_1,t_2}$ {\it and (2.9), (2.10) for the}
${\cal J}${\it - and} ${\cal I}${\it -functions, respectively. }

{\bf REMARK 2.2.} We will see later that the class numbers $h\left
(t_1^2-4,t_2^2-4,f\right)$ are
finite, see Lemma 3.1 and Remark 3.1. The sum (2.14) is
actually finite, because for large enough $\left|f\right|$ the
${\cal I}-$function there is $0$. This can be seen easily from (2.10)
and (2.11).

For the proof of Lemma 2.2 we need a few preliminary
lemmas. To state the first one, we give some definitions
using the notations of Lemma 2.2.

Write $G:=\Gamma_{t_1}\times\Gamma_{t_2},$ and let $G_0$ be the set of those elements $\left
(\gamma_1,\gamma_2\right)\in G$ for
which the set of fixed points on ${\bf R}$ of $\gamma_1$ and of $
\gamma_2$ are the same. If $\left(\gamma_1,\gamma_2\right),\left(
\gamma_1^{\ast},\gamma_2^{\ast}\right)\in G$, we say that $\left(
\gamma_1,\gamma_2\right)$ and
$\left(\gamma_1^{\ast},\gamma_2^{\ast}\right)$ are $SL_2({\bf Z})$-equivalent if there is an element $
\tau\in SL_2({\bf Z})$ such that
$\tau^{-1}\gamma_i\tau =\gamma_i^{\ast}$ for $i=1,2$. We denote by $
G_0^{\ast}$ a complete set of representatives of the $SL_2({\bf Z}
)$-equivalence classes of $G_0$, and by $\left(G\setminus G_0\right
)^{\ast}$ a complete set of representatives of the
$SL_2({\bf Z})$-equivalence classes of $G\setminus G_0$.

{\bf LEMMA 2.3.} {\it Let} $t_1,$$t_2>2$ {\it be integers and let} $
m_1${\it ,} $m_2$ {\it be compactly supported bounded functions on} $
[0,\infty )${\it . Recall the notation (2.8). We have that}
$$\int_{{\cal F}}M_{t_1,m_1}(z)M_{t_2,m_2}(z)d\mu_z\eqno (2.15)$$
{\it equals the sum of }
$$\sum_{\left(\gamma_1,\gamma_2\right)\in G_0^{\ast}}\int_{C\left
(\gamma_1\right)\setminus\bbb H}m_1\left(z,\gamma_1z\right)m_2\left
(z,\gamma_2z\right)d\mu_z$$
{\it and}
$$\sum_{\left(\gamma_1,\gamma_2\right)\in\left(G\setminus G_0\right
)^{\ast}}\int_{\bbb H}m_1\left(z,\gamma_1z\right)m_2\left(z,\gamma_
2z\right)d\mu_z.$$
{\bf REMARK 2.3.} To avoid confusion we emphasize that $G\setminus
G_0$
denotes set difference, while $C\left(\gamma_1\right)\setminus\bbb
H$ denotes
quotient on the left.

{\it Proof.\/} An element $\gamma\in SL_2({\bf Z})$ with $\hbox{\rm tr$\left
(\gamma\right)$}>2$ determines a
hyperbolic transformation of $\bbb H$, see Section 1.5 of [I]. Hence $
\gamma$ has two different fixed points on ${\bf R}$. Assume that $\gamma_1\in\Gamma_{t_
1}$, $\gamma_2\in\Gamma_{t_2}$, $\tau\in SL_2({\bf Z})$ and
$$\tau^{-1}\gamma_1\tau =\gamma_1,\quad\tau^{-1}\gamma_2\tau =\gamma_
2.\eqno (2.16)$$
It is clear by (2.4) that if $\left(\gamma_1,\gamma_2\right)\in G
\setminus G_0$, then (2.16) is true
if and only if $\tau =\pm\left(\matrix{1&0\cr
0&1\cr}
\right)$. If $\left(\gamma_1,\gamma_2\right)\in G_0$, then by (2.4)
we see that $C\left(\gamma_1\right)=C\left(\gamma_2\right)$, and (2.16) is true if and only
if $\tau\in C\left(\gamma_1\right)$.

By the definitions we see that (2.15) equals
$$\sum_{\gamma_1\in\Gamma_{t_1}}\sum_{\gamma_2\in\Gamma_{t_2}}\int_
{\cal F}m_1\left(z,\gamma_1z\right)m_2\left(z,\gamma_2z\right)d\mu_z.$$
We partition $G$ into $SL_2({\bf Z})$-equivalence classes. Since for
$\tau\in SL_2({\bf Z})$ we have that
$$\int_{\cal F}m_1\left(z,\tau^{-1}\gamma_1\tau z\right)m_2\left(z,\tau^{
-1}\gamma_2\tau z\right)d\mu_z=\int_{\tau {\cal F}}m_1\left(z,\gamma_1z\right
)m_2\left(z,\gamma_2z\right)d\mu_z,$$
our considerations above give the lemma.

{\bf LEMMA 2.4.} {\it Let} $\gamma_1=\left(\matrix{a&b\cr
c&d\cr}
\right),$ $\gamma_2=\left(\matrix{A&B\cr
C&D\cr}
\right)$ {\it be hyperbolic elements of} $SL_2({\bf R})$, {\it assume that the set of fixed points of} $
\gamma_1$ {\it and the set of fixed points of} $\gamma_2$ {\it are disjoint. Let}
$$F:=F\left(\gamma_1,\gamma_2\right)={{\left(d-a\right)\left(D-A\right
)+2bC+2Bc}\over {\sqrt {\left(d+a\right)^2-4}\sqrt {\left(D+A\right
)^2-4}}}.\eqno (2.17)$$
{\it Let us write} $t_1=a+d$, $t_2=A+D${\it , and assume} $t_1,t_
2>2$. {\it Let }
$m_1${\it ,} $m_2$ {\it be compactly supported bounded functions on} $
[0,\infty )$ {\it and use the notation (2.7). Then we have that }
$$\int_{\bbb H}m_1\left(z,\gamma_1z\right)m_2\left(z,\gamma_2z\right
)d\mu_z={\cal I}\left(t_1,t_2,F,m_1,m_2\right),\eqno (2.18)$$
{\it where the} ${\cal I}${\it -function is defined in (2.10) and (2.11). }

{\it Proof.\/} First note that it is easy to check that
$F\left(\gamma_1,\gamma_2\right)=F\left(\tau^{-1}\gamma_1\tau ,\tau^{
-1}\gamma_2\tau\right)$ for $\tau\in SL_2({\bf R})$. Since (2.18) also remains the same if we write $
\tau^{-1}\gamma_1\tau$ and
$\tau^{-1}\gamma_2\tau$ in place of $\gamma_1$ and $\gamma_2$, and since we can choose $
\tau$
in such a way that $\tau^{-1}\gamma_1\tau$ is diagonal, for the
proof of the lemma we may assume that $\gamma_1$ is diagonal.

So assume that $b=c=0$. Then by the conditions we have
$BC\neq 0$. It can be easily computed by the definitions that
$$u\left(z,\gamma_1z\right)={{\left(d-a\right)^2\left|z\right|^2}\over {
4\hbox{\rm $\hbox{\rm Im}^2$$z$}}},\;u\left(z,\gamma_2z\right)={{\left
|Cz^2+\left(D-A\right)z-B\right|^2}\over {4\hbox{\rm $\hbox{\rm Im}^
2$$z$}}}.$$
Hence if $z=x+iy$, then using $ad=1$ and $AD-BC=1$ we get by
some computations that
$$u\left(z,\gamma_1z\right)={{\left(a+d\right)^2-4}\over 4}+{{\left
(d-a\right)^2x^2}\over {4\hbox{\rm $y^2$}}},\eqno (2.19)$$
$$u\left(z,\gamma_2z\right)={{\left(A+D\right)^2-4}\over 4}+{{\left
(Cx^2+\left(D-A\right)x-B+Cy^2\right)^2}\over {4\hbox{\rm $y^2$}}}
.\eqno (2.20)$$
Since $m(z,w)$ is defined through the function $u$, we will
be able to compute the left-hand side of (2.18) using (2.19)
and (2.20).

Let us use the substitution
$$q:={x\over {\hbox{\rm $y$}}},\;r:={{Cx^2+\left(D-A\right)x-B+Cy^
2}\over {\hbox{\rm $y$}}}.\eqno (2.21)$$
Then the determinant of the Jacobi matrix ${{dqdr}\over {dxdy}}$ is
$$\det\left(\matrix{{1\over y}&{{2Cx+\left(D-A\right)}\over {\hbox{\rm $
y$}}}\cr
-{x\over {\hbox{\rm $y$}^2}}&C-{{Cx^2+\left(D-A\right)x-B}\over {\hbox{\rm $
y$}^2}}\cr}
\right)={{B+Cx^2+Cy^2}\over {\hbox{\rm $y$}^3}}.$$
It is not hard to check that
$${{B+Cx^2+Cy^2}\over {\hbox{\rm $y$}}}=r-\left(D-A\right)q+{{2B}\over
y}\eqno (2.22)$$
and
$${B\over {\hbox{\rm $y$}^2}}+{{r-\left(D-A\right)q}\over y}-\left
(C+Cq^2\right)=0.\eqno (2.23)$$
From (2.22) and (2.23) we get
$$\left({{B+Cx^2+Cy^2}\over {\hbox{\rm $y$}}}\right)^2=\left(\left
(D-A\right)q-r\right)^2+4B\left(C+Cq^2\right).$$
Hence if we want to compute the left-hand side of (2.18) by the substitution
(2.21), then on the one hand we see that for $q$ and
$r$ we have the condition
$$\left(\left(D-A\right)q-r\right)^2+4B\left(C+Cq^2\right)\ge 0.\eqno
(2.24)$$
On the other hand, in the case $BC>0$ we see from the quadratic
equation (2.23) that for every real $q$ and $r$ satisfying (2.24) there is exactly one
$y>0$ and real $x$ satisfying (2.21). Similarly, in the
case $BC<0$ we see from the quadratic equation (2.23) that
we must have ${{\left(D-A\right)q-r}\over B}>0$, i.e. combined with (2.24) we
must have
$${{\left(D-A\right)q-r}\over B}\ge 2\sqrt {-{{C+Cq^2}\over B}}.\eqno
(2.25)$$
If the left-hand side of (2.25) is larger than the
right-hand side, then we have two positive solutions
of (2.23) in ${1\over y}$. If (2.25) holds with equality, then we
have a double root.

Putting everything together we see that the left-hand side
of (2.18) equals
$$\int_{-\infty}^{\infty}\int_{A_q}f\left(r,q\right)drdq\eqno (2.
26)$$
for $BC>0$, and the left-hand side of (2.18) equals
$$2\int_{-\infty}^{\infty}\int_{A_q^{+}}f\left(r,q\right)drdq\eqno
(2.27)$$
for $BC<0$, where

$$f\left(r,q\right):={{m_1\left({{\left(a+d\right)^2-4}\over 4}+{{\left
(d-a\right)^2q^2}\over 4}\right)m_2\left({{\left(A+D\right)^2-4}\over
4}+{{r^2}\over 4}\right)}\over {\sqrt {\left(\left(D-A\right)q-r\right
)^2+4B\left(C+Cq^2\right)}}},$$
$$A_q:=\left\{r\in {\bf R}:\;\left(\left(D-A\right)q-r\right)^2+4
B\left(C+Cq^2\right)\ge 0\right\},$$
and for $BC<0$ we write
$$A_q^{+}:=\left\{r\in {\bf R}:\;{{\left(D-A\right)q-r}\over B}\ge
2\sqrt {-{{C+Cq^2}\over B}}\right\}.$$
For $BC<0$ we define also
$$A_q^{-}:=\left\{r\in {\bf R}:\;{{\left(D-A\right)q-r}\over B}\le
-2\sqrt {-{{C+Cq^2}\over B}}\right\}.$$
We see that for $BC<0$ we have that $-A_{-q}^{+}=A_q^{-}$, so, since
$f\left(r,q\right)=f\left(-r,-q\right)$, for $BC<0$ we have that (2.27) equals
$$\int_{-\infty}^{\infty}\int_{A_q^{+}}f\left(r,q\right)drdq+\int_{
-\infty}^{\infty}\int_{A_q^{-}}f\left(r,q\right)drdq.$$
Since $A_q$ is the disjoint union of $A_q^{+}$ and $A_q^{-}$ for $
BC<0$, we finally get that the left-hand
side of (2.18) equals (2.26) also in the case $BC<0$. Apply the substitution
$$S=q,\;T={{-r}\over {\sqrt {\left(D+A\right)^2-4}}}$$
in (2.26). Recalling (2.17), $b=c=0$ and $BC\neq 0$ we have that
$$\left(a+d\right)^2-4=\left(a-d\right)^2,\quad\left(D-A\right)^2
+4BC=\left(D+A\right)^2-4,\quad F^2={{\left(D-A\right)^2}\over {\left
(D+A\right)^2-4}}\neq 1.$$
Taking into account the fact that ${\cal I}\left(t_1,t_2,F,m_1,m_
2\right)$ is even
in $F$, we get equation (2.18) for the case $b=c=0$. But we
have seen that then the lemma is completely proved.

{\bf LEMMA 2.5.} {\it Let} $\gamma_1\in\Gamma_{t_1},$ $\gamma_2\in
\Gamma_{t_2},$ {\it where} $t_i>2$ {\it for} $i=1,2$. {\it Assume that the set of fixed points of} $
\gamma_1$ {\it and the set of fixed points of} $\gamma_2$ {\it are the same. Then we have that }
$$\int_{C\left(\gamma_1\right)\setminus\bbb H}m_1\left(z,\gamma_1
z\right)m_2\left(z,\gamma_2z\right)d\mu_z={\cal J}\left(t_1,t_2,m_
1,m_2\right)\left|\log N\left(\gamma_0\right)\right|,\eqno (2.28)$$
{\it where} $\gamma_0\in SL_2({\bf Z})$ {\it is a generator of the centralizer} $
C\left(\gamma_1\right)${\it , see (2.5). The} ${\cal J}${\it -function is defined in (2.9).}

{\it Proof.\/} We may assume that $N\left(\gamma_0\right)>1$. We can choose $
\tau\in SL_2({\bf R})$
in such a way that $\tau^{-1}\gamma_i\tau z=\lambda_iz$ for every $
z\in\bbb H$ and
$0\le i\le 2$ with $\lambda_0=N\left(\gamma_0\right)$ and $\lambda_
i=N\left(\gamma_i\right)^{\epsilon_i}$ for $i=1,2$,
where $\epsilon_i\in\left\{-1,1\right\}$. The fundamental domain of the group
$\tau^{-1}C\left(\gamma_1\right)\tau$ in $\bbb H$ is the subset $\left
\{\,1\le\left|z\right|<N\left(\gamma_0\right)\right\}$. Then we have that the left-hand
side of (2.28) equals
$$\int_{\left\{z\in {\bbb H}:\,1\le\left|z\right|<N\left(\gamma_0\right)\right
\}}m_1\left(z,\tau^{-1}\gamma_1\tau z\right)m_2\left(z,\tau^{-1}\gamma_
2\tau z\right)d\mu_z.\eqno (2.29)$$
By the substitution $z=re^{i\left({{\pi}\over 2}+\theta\right)}$ with $
1\le r<N\left(\gamma_0\right)$,
$-{{\pi}\over 2}<\theta <{{\pi}\over 2}$ we get using $d\mu_z={{d
rd\theta}\over {r\cos^2\theta}}$ that $(2.29)$ equals
$$\int_{-\pi /2}^{\pi /2}\int_1^{N\left(\gamma_0\right)}m_1\left({{
\lambda_1+\lambda_1^{-1}-2}\over {4\cos^2\theta}}\right)m_2\left({{
\lambda_2+\lambda_2^{-1}-2}\over {4\cos^2\theta}}\right){{drd\theta}\over {
r\cos^2\theta}}.$$
It is clear that $\lambda_1^{1/2}+\lambda_1^{-1/2}=t_1$, $\lambda_
2^{1/2}+\lambda_2^{-1/2}=t_2$, since the trace is invariant under conjugation. The lemma is proved.

{\it Proof of Lemma 2.2.\/} We use the one-to-one correspondence betweeen
$\Gamma_{t_i}$ and $\scr{Q}_{t_i^2-4}$ for $i=1,2$.

Using the notations of Lemma 2.3, we see that if $\left(\gamma_1,
\gamma_2\right)\in G$, then $\left(\gamma_1,\gamma_2\right)\in G_
0$ holds if
and only if $Q_{\gamma_1}=\lambda Q_{\gamma_2}$ with some $\lambda
\in {\bf Q}$. Indeed,
$Q_{\gamma_1}=\lambda Q_{\gamma_2}$ holds with some $\lambda\in {\bf Q}$ if and only if the
polynomials $Q_{\gamma_1}(X,1)$ and $Q_{\gamma_2}(X,1)$ have the same roots,
i.e. if and only if $\gamma_1$ and $\gamma_2$ have the same fixed points.

It is clear that $\left(\gamma_1,\gamma_2\right),\left(\gamma_1^{
\ast},\gamma_2^{\ast}\right)\in G$ are $SL_2({\bf Z})$-equivalent if and only if there is a $
\tau\in SL_2({\bf Z})$
such that $\left(Q_{\gamma_1}^{\tau},Q_{\gamma_1}^{\tau}\right)=\left
(Q_{\gamma_1^{\ast}},Q_{,\gamma_2^{\ast}}\right)$. We note also that if $\left
(\gamma_1,\gamma_2\right)\in G\setminus G_0$, then for the quantity
$F\left(\gamma_1,\gamma_2\right)$ defined in (2.17) we have $F\left
(\gamma_1,\gamma_2\right)={f\over {\sqrt {t_1^2-4}\sqrt {t_2^2-4}}}$
with the notations
$$Q_{\gamma_1}\left(X,Y\right)=A_1X^2+B_1XY+C_1Y^2,\;Q_{\gamma_2}\left
(X,Y\right)=A_2X^2+B_2XY+C_2Y^2,\eqno (2.30)$$
$$f=B_1B_2-2A_1C_2-2A_2C_1.\eqno (2.31)$$
We show that if $\left(\gamma_1,\gamma_2\right)\in G\setminus G_0$, then
$f^2\neq\left(t_1^2-4\right)\left(t_2^2-4\right)$. Indeed, writing $
d_i:=t_i^2-4$ and using
$d_i=B_i^2-4A_iC_i$ for $i=1,2$, we easily get from (2.30) and
(2.31) that
$$d_2^2A_1^2-2fd_2A_1A_2+d_1d_2A_2^2=d_2\left(B_2A_1-B_1A_2\right
)^2,$$
$$d_2^2C_1^2-2fd_2C_1C_2+d_1d_2C_2^2=d_2\left(B_2C_1-B_1C_2\right
)^2.$$
Assume $f^2=d_1d_2$. Then the left-hand sides above are squares, and
since $d_2=t_2^2-4$ cannot be a square, we get
$B_2A_1-B_1A_2=0$, $B_2C_1-B_1C_2=0$. One has the identity
$$\left(A_1C_2-A_2C_1\right)^2-\left(A_1B_2-A_2B_1\right)\left(B_
1C_2-B_2C_1\right)={{f^2-\prod_{i=1}^2\left(B_i^2-4A_iC_i\right)}\over
4}\eqno (2.32)$$
with $f$ defined in (2.31). We get also $A_1C_2-A_2C_1=0$ from
(2.32). It follows that the vectors $\left(A_1,B_1,C_1\right)$ and
$\left(A_2,B_2,C_2\right)$ are linearly dependent, hence $\left(\gamma_
1,\gamma_2\right)\in G_0$, which is a contradiction.

We note finally that if $\left(\gamma_1,\gamma_2\right)\in G\setminus
G_0$, then the set of fixed points of $\gamma_1$ and the set of fixed points of $
\gamma_2$ are disjoint.

By these considerations, applying Lemmas 2.3, 2.4, 2.5 and the
definitions we obtain the lemma.

{\bf 3. Estimates on the number of equivalence classes of quadratic forms.}
\medskip

Recall the definition of ${\cal Q}_{d_1,d_2,t}$, $h\left(d_1,d_2,
t\right)$, see (2.1), (2.2),
and the paragraph below (2.2). In this section we will give several upper bounds for
$h\left(d_1,d_2,t\right)$ itself and for certain sums containing $
h\left(d_1,d_2,t\right)$.

Let ${\rm g}{\rm c}{\rm d}\left(n_1,n_2,\ldots ,n_r\right)$ be the greatest common divisor of
the integers $n_1,n_2,\ldots ,n_r$. The integer part of a real number $
x$ is denoted by $[x]$.

{\bf 3.1. A general upper bound for} $h\left(d_1,d_2,t\right)${\bf .} Our aim in this
subsection is to prove Lemma 3.1. The upper bound we give for $h\left
(d_1,d_2,t\right)$ will be smaller than
$\left(1+\left|d_1d_2t\right|\right)^{\epsilon}$ for any fixed $\epsilon
>0$ in many cases. This is not always true, but the exceptions are rare.

For any finite set of integers $n_1,n_2,\ldots ,n_r$ we write
$$S\left(n_1,n_2,\ldots ,n_r\right)=\max\left\{k\ge 1:\;k^2|{\rm g}
{\rm c}{\rm d}\left(n_1,n_2,\ldots ,n_r\right)\right\}.\eqno (3.1
)$$
Denote by $\tau\left(n\right)$ the number of divisors and by $\omega\left
(n\right)$ the
number of distinct prime divisors of a nonzero integer $n$.

{\bf LEMMA 3.1.} {\it Assume that} $d_1,d_2,t\in {\bf Z}$, {\it and} $
d_i$ {\it is not a square of an integer} $(i=1,2)$. {\it Assume also that} $
t^2-d_1d_2\neq 0$. {\it Then we have that}
$$h\left(d_1,d_2,t\right)\le C2^{\omega\left(t^2-d_1d_2\right)}\tau\left(t^2-d_1d_2\right)S\left(d_1,d_2,t^2\right),\eqno (3.2)$$
{\it where} $C>0$ {\it is an absolute constant.}

{\bf REMARK 3.1.} Even the finiteness of $h\left(d_1,d_2,t\right)$ is not
completely obvious. For a short proof of this fact using the
theory of algebraic groups see Appendix I of [M].

Note that for the case $d_i=t_i^2-4$ for $i=1,2$ with
integers $t_i>2$, which is our primary interest, one can
give a trivial upper bound for $h\left(d_1,d_2,t\right)$ using Lemma 2.2. For simplicity let us consider the case when
$X\le t_1^2-4,t_2^2-4,t\le 2X$ with a large real number $X$. Then
choosing $m_1$ and also $m_2$ in Lemma 2.2 to be the characteristic function of the
interval $\left[0,CX\right]$ with a suitable absolute constant $C$ one
can show the trivial bound $h\left(t_1^2-4,t_2^2-4,t\right)\ll X$. Indeed, the coefficient of $
h\left(t_1^2-4,t_2^2-4,t\right)$ is bounded from
below by a positive constant, every term in (2.14) and
(2.13) is nonnegative, and one can prove that (2.12) equals
$O\left(X\right)$ in this case. The estimate (3.2) gives better than $
h\left(t_1^2-4,t_2^2-4,t\right)\ll X$ even in the worst case, when
$S\left(t_1^2-4,t_2^2-4,t^2\right)$ is as large as $\sqrt X$. But the $
S$-function
is often much smaller than $\sqrt X$, so the bound given in
Lemma 3.1 is much stronger than the trivial bound.

To prepare the proof of Lemma 3.1 we need two
preliminary lemmas. We introduce the notation
$$C_{d_1,d_2,t}:=\left\{\left(x,y\right)\in {\bf R}^2:\;d_2x^2+d_
1y^2-2txy=1\right\}.$$
In the first lemma we prove general statements for any
two different points of $C_{d_1,d_2,t}$. In the second one we
show that if we have any element of ${\cal Q}_{d_1,d_2,t}$, then we
can parametrize the rational points of $C_{d_1,d_2,t}$.

{\bf LEMMA 3.2.} {\it Let} $d_1,d_2,t$ {\it be as in Lemma 3.1, assume that }
$\left(x_i,y_i\right)\in C_{d_1,d_2,t}$ {\it for} $i=1,2$ {\it and} $\left
(x_1,y_1\right)\neq\left(x_2,y_2\right)${\it . Then we have}
$$d_2\left(x_1-x_2\right)^2+d_1\left(y_1-y_2\right)^2-2t\left(x_1
-x_2\right)\left(y_1-y_2\right)\neq 0\eqno (3.3)$$
{\it and}
$$\left(d_1y_1-tx_1\right)\left(y_1-y_2\right)+\left(d_2x_1-ty_1\right
)\left(x_1-x_2\right)\neq 0.\eqno (3.4)$$
{\it Proof.\/} Let $S_1$ and $S_2$ be the quantities appearing in (3.3)
and (3.4), respectively. One can check the identities
$$S_2=\left(\matrix{y_1&x_1\cr}
\right)\left(\matrix{d_1&-t\cr
-t&d_2\cr}
\right)\left(\matrix{y_1-y_2\cr
x_1-x_2\cr}
\right)\eqno (3.5)$$
and
$$2S_2+\sum_{i=1}^2\left(-1\right)^i\left(d_2x_i^2+d_1y_i^2-2tx_i
y_i\right)=S_1.$$
Since $\left(x_i,y_i\right)\in C_{d_1,d_2,t}$ for $i=1,2$, we have $
S_1=2S_2$. Hence it
is enough to show that $S_1\neq 0$. Assume for a contradiction
that $S_1=0$. Then the right-hand side of (3.5) is $0$, but this is true also by
exchanging the role of $\left(x_1,y_1\right)$ and $\left(x_2,y_2\right
)$, so we get
$$\left(\matrix{y_1&x_1\cr
y_2&x_2\cr}
\right)\left(\matrix{d_1&-t\cr
-t&d_2\cr}
\right)\left(\matrix{y_1-y_2\cr
x_1-x_2\cr}
\right)=\left(\matrix{0\cr
0\cr}
\right).$$
The vector $\left(\matrix{y_1-y_2\cr
x_1-x_2\cr}
\right)$ is nonzero and $\det\left(\matrix{d_1&-t\cr
-t&d_2\cr}
\right)\neq 0$
by $t^2-d_1d_2\neq 0$, so we must have $\det\left(\matrix{y_1&x_1\cr
y_2&x_2\cr}
\right)=0$. Hence
$\left(x_2,y_2\right)=\lambda\left(x_1,y_1\right)$ with some constant $
\lambda\neq 1$, so $S_1=\left(1-\lambda\right)^2\neq 0$ by our
assumptions. This is a contradiction, the lemma is proved.

{\bf LEMMA 3.3.} {\it Let} $d_1,d_2,t$ {\it be as in Lemma 3.1, and let }
$Q_i\left(X,Y\right)=A_iX^2+B_iXY+C_iY^2$ $(i=1,2)$ {\it be such that} $\left
(Q_1,Q_2\right)\in {\cal Q}_{d_1,d_2,t}$. {\it Assume that} $A_1B_2-A_2B_
1\neq 0$. {\it Define}
$$R\left(X,Y\right)=\left(A_1B_2-A_2B_1\right)X^2+2XY\left(A_1C_2
-A_2C_1\right)+Y^2\left(B_1C_2-B_2C_1\right).$$
{\it Then for} $x,y\in {\bf Q}$ {\it the following two statements are equivalent.}

{\it (i) We have} $\left(x,y\right)\in C_{d_1,d_2,t}$.

{\it (ii) There are} $a,b\in {\bf Q}$ {\it such that} $R\left(a,b\right
)\neq 0$ {\it and writing }
$x_{a,b}:={{Q_1\left(a,b\right)}\over {R\left(a,b\right)}}$, $y_{
a,b}:={{Q_2\left(a,b\right)}\over {R\left(a,b\right)}}$ {\it we have} $
x=x_{a,b}${\it ,} $y=y_{a,b}${\it . }

{\it Proof.\/} By straightforward computation using the definitions we get
the identity
$$d_2\left(Q_1\left(a,b\right)\right)^2+d_1\left(Q_2\left(a,b\right
)\right)^2-2tQ_1\left(a,b\right)Q_2\left(a,b\right)=\left(R\left(
a,b\right)\right)^2.\eqno (3.6)$$
Introduce the abbreviations
$$a_1={{A_1}\over {A_1B_2-A_2B_1}},\;a_2={{A_2}\over {A_1B_2-A_2B_
1}}.\eqno (3.7)$$
Note that writing $a=1$, $b=0$ in (3.6) we get
$\left(a_1,a_2\right)\in C_{d_1,d_2,t}$. We first assume (ii). Then (i) follows at
once from (3.6).

We now assume (i). If $\left(x,y\right)=\left(a_1,a_2\right)$, then we can
take $a=1$, $b=0$. So let us assume that $\left(x,y\right)\neq\left
(a_1,a_2\right)$. It is easy to see that if $a,b\in {\bf Q}$, then
$$Q_1\left(a,b\right)-a_1R\left(a,b\right)={{b\left(a\alpha +b\beta\right
)}\over {\left(A_1B_2-A_2B_1\right)}},\quad Q_2\left(a,b\right)-a_
2R\left(a,b\right)={{b\left(\gamma a+\delta b\right)}\over {\left
(A_1B_2-A_2B_1\right)}}\eqno (3.8)$$
with
$$\alpha :=B_1\left(A_1B_2-A_2B_1\right)+2A_1\left(A_2C_1-A_1C_2\right
)=tA_1-d_1A_2,\eqno (3.9)$$
$$\beta :=C_1\left(A_1B_2-A_2B_1\right)+A_1\left(C_1B_2-C_2B_1\right
),\eqno (3.10)$$
$$\gamma :=B_2\left(A_1B_2-A_2B_1\right)+2A_2\left(A_2C_1-A_1C_2\right
)=-tA_2+d_2A_1,\eqno (3.11)$$
$$\delta :=C_2\left(A_1B_2-A_2B_1\right)+A_2\left(C_1B_2-C_2B_1\right
).\eqno (3.12)$$
Let us write $g:=\left(\matrix{\alpha&\beta\cr
\gamma&\delta\cr}
\right)$. Now, one can compute that
$$\det g=2\left(A_2B_1-A_1B_2\right)\left(\left(A_1C_2-A_2C_1\right
)^2-\left(A_1B_2-A_2B_1\right)\left(B_1C_2-B_2C_1\right)\right).$$
The last bracket equals ${{t^2-d_1d_2}\over 4}$ by (2.32) and (2.31). Hence
$t^2-d_1d_2\neq 0$ and $A_1B_2-A_2B_1\neq 0$ imply $\det g\neq 0$.

Let us take $a,b\in {\bf Q}$ in the following way:
$$\left(\matrix{a\cr
b\cr}
\right)=\left(\matrix{\delta&-\beta\cr
-\gamma&\alpha\cr}
\right)\left(\matrix{x-a_1\cr
y-a_2\cr}
\right).\eqno (3.13)$$
By (3.8) and (3.13) we then easily get
$$\left(\matrix{Q_1\left(a,b\right)-a_1R\left(a,b\right)\cr
Q_2\left(a,b\right)-a_2R\left(a,b\right)\cr}
\right)={{b\det g}\over {A_1B_2-A_2B_1}}\left(\matrix{x-a_1\cr
y-a_2\cr}
\right).\eqno (3.14)$$
So this is true if $\left(x,y\right)\neq\left(a_1,a_2\right)$, and $
a$, $b$ are defined
by (3.9)-(3.13).

Assume that $b=0$. Then by (3.9), (3.11) and (3.13) we get
$$\left(d_1A_2-tA_1\right)\left(y-a_2\right)+\left(-tA_2+d_2A_1\right
)\left(x-a_1\right)=0.$$
By (3.7), $\left(a_1,a_2\right),\left(x,y\right)\in C_{d_1,d_2,t}$ and $\left
(x,y\right)\neq\left(a_1,a_2\right)$ this
contradicts (3.4). So we have $b\neq 0$.

Assume that $R\left(a,b\right)=0$. Then (3.6) and (3.14) imply that
$$d_2\left(x-a_1\right)^2+d_1\left(y-a_2\right)^2-2t\left(x-a_1\right
)\left(y-a_2\right)=0.$$
But this contradicts (3.3). So we have $R\left(a,b\right)\neq 0$.

Then (3.14) clearly implies
$$\left(\matrix{x_{a,b}-a_1\cr
y_{a,b}-a_2\cr}
\right)={{b\det g}\over {R\left(a,b\right)\left(A_1B_2-A_2B_1\right
)}}\left(\matrix{x-a_1\cr
y-a_2\cr}
\right).$$
Hence we have $\left(x_{a,b},y_{a,b}\right)\neq\left(a_1,a_2\right
)$, since we
assumed $\left(x,y\right)\neq\left(a_1,a_2\right)$.  We would like to show that
$\left(x,y\right)=\left(x_{a,b},y_{a,b}\right)$. If this is false, then $\left
(a_1,a_2\right)$,
$\left(x_{a,b},y_{a,b}\right)$ and $\left(x,y\right)$ are three pairwise different
points lying on a line and all of these three points belong
to $C_{d_1,d_2,t}$. Hence we have that the equation
$$d_2\left(a_1+q\left(x-a_1\right)\right)^2+d_1\left(a_2+q\left(y
-a_2\right)\right)^2-2t\left(a_1+q\left(x-a_1\right)\right)\left(
a_2+q\left(y-a_2\right)\right)=1$$
has three different real solutions in $q$. The coefficient of
$q^2$ is nonzero by (3.3), so this is a contradiction. The
lemma is proved.

For the next proof we need the following notation. If $p$ is a prime and $
n\neq 0$ is an
integer, let us denote by $\nu_p\left(n\right)$ the largest nonnegative
integer such that $p^{\nu_p\left(n\right)}$ divides $n.$

{\it Proof of Lemma 3.1.\/} If ${\cal Q}_{d_1,d_2,t}$ is empty, then
$h\left(d_1,d_2,t\right)=0$ and there is nothing to prove. So assume in the sequel that $
{\cal Q}_{d_1,d_2,t}\neq\emptyset$. We divide the proof into
four parts, we formulate them as claims.

{\bf CLAIM A.} Let $Q_i\left(X,Y\right)=A_iX^2+B_iXY+C_iY^2$ $(i=
1,2)$ be two quadratic forms such that
$\left(Q_1,Q_2\right)\in {\cal Q}_{d_1,d_2,t}$. Then replacing $\left
(Q_1,Q_2\right)$ by an element in its
$SL_2({\bf Z})$-equivalence class we can achieve $B_2A_1-B_1A_2\neq
0$.

{\it Proof of Claim A.\/} We first show the weaker statement that replacing $\left
(Q_1,Q_2\right)$ by an element in its $SL_2({\bf Z})$-equivalence class we may assume that
$\left(B_1,B_2\right)\neq\left(0,0\right)$. If $\tau =\left(\matrix{
1&b\cr
0&1\cr}
\right)$, then we have  $ $
$$Q_i^{\tau}\left(X,Y\right)=Q_i\left(X+bY,Y\right)=A_iX^2+\left(
B_i+2A_ib\right)XY+C_i^{\ast}Y^2.\eqno (3.15)$$
for $i=1,2$ with some $C_i^{\ast}$. If $\left(B_1+2A_1b,B_2+2A_2b\right
)=\left(0,0\right)$
for every integer $b$, then $A_i=B_i=0$ for $i=1,2$. But this
is impossible, since this would imply $d_1=d_2=t=0$, but this contradicts $
t^2-d_1d_2\neq 0$.

Hence we may assume that $\left(B_1,B_2\right)\neq\left(0,0\right
)$. Let $B_1\neq 0$,
say. Assume for a contradiction that $B_2A_1-B_1A_2=0$ and
$B_2C_1-B_1C_2=0$. Then $\left(A_2,C_2\right)=\lambda\left(A_1,C_
1\right)$
with $\lambda =B_2/B_1$, hence $C_2A_1-C_1A_2=0$. So the matrix
$\left(\matrix{A_1&B_1&C_1\cr
A_2&B_2&C_2\cr}
\right)$ has rank $1$, hence its lines are linearly
dependent. But this contradicts $t^2-d_1d_2\neq 0$.

So we may assume that $B_2A_1-B_1A_2\neq 0$ or
$B_2C_1-B_1C_2\neq 0$. But applying the matrix $\left(\matrix{0&-
1\cr
1&0\cr}
\right)\in SL_2({\bf Z})$
we can exchange the roles of $A_i$ and $C_i$. Claim A follows.

{\bf CLAIM B.} We have $C_{d_1,d_2,t}\neq\emptyset$.

{\it Proof of Claim B.\/} By Claim A there is an element $\left(Q_
1,Q_2\right)\in {\cal Q}_{d_1,d_2,t}$ such that for their coefficients we have
$B_2A_1-B_1A_2\neq 0$. Taking $a=1$, $b=0$ in Lemma 3.3 (ii) we see
by that lemma that there are numbers $x,y\in {\bf Q}$ such that
$\left(x,y\right)\in C_{d_1,d_2,t}$. Claim B is proved.

{\bf CLAIM C.} There exists a subset ${\cal A}$ of ${\bf Z}^2$ of size $
2\tau\left(t^2-d_1d_2\right)$ such that in every
$SL_2({\bf Z})$-equivalence class of ${\cal Q}_{d_1,d_2,t}$ there is an element
$\left(Q_1,Q_2\right)$ with coefficients  $Q_i\left(X,Y\right)=A_
iX^2+B_iXY+C_iY^2$
$(i=1,2)$ such that $\left(A_1,A_2\right)\in {\cal A}$.

{\it Proof of Claim C.\/} Fix $x,y\in {\bf Q}$ such that $\left(x
,y\right)\in C_{d_1,d_2,t}$,
this is possible by Claim B. We fix also integers $m$ and $n$
satisfying ${\rm g}{\rm c}{\rm d}\left(m,n\right)=1$ and an $s\in
{\bf Q}$, $s\neq 0$ such that
$x=sm$, $y=sn$. Let us take an arbitrary $SL_2({\bf Z})$-equivalence class of
${\cal Q}_{d_1,d_2,t}$. We know by Claim A that we can take an element $\left
(Q_1,Q_2\right)$ in this equivalence
class such that we have $B_2A_1-B_1A_2\neq 0$ for their
coefficients. Then it follows from Lemma 3.3 that there
are $a,b\in {\bf Q}$ such that $Q_1\left(a,b\right)=qx$, $Q_2\left
(a,b\right)=qy$ with
some $q\in {\bf Q}$, $q\neq 0$. We may clearly assume here that
$a,b\in {\bf Z}$ and $\left(a,b\right)=1$. Taking $\left(\matrix{
a&c\cr
b&d\cr}
\right)\in SL_2({\bf Z})$
with some suitable $c$ and $d$ we then see that replacing $\left(
Q_1,Q_2\right)$ by an element in
its $SL_2({\bf Z})$-equivalence class we may assume that for
their coefficients we have $A_1=rx$, $A_2=ry$ with some $r\in {\bf Q}$, $
r\neq 0$. Observe also that
$${\rm g}{\rm c}{\rm d}\left(A_1,A_2\right)\left|t^2-d_1d_2\right
.\eqno (3.16)$$
follows from the definition of $d_1$, $d_2$ and $t$. Then recalling $
x=sm$, $y=sn$ we see that $A_1=\delta m$, $A_2=\delta n$ with some
integer $\delta$ dividing $t^2-d_1d_2$. Claim C is proved.

{\bf CLAIM D.} Let the integers $A_1$ and $A_2$ be given and assume that there are $
N$
inequivalent elements in ${\cal Q}_{d_1,d_2,t}$ of the form $Q_i\left
(X,Y\right)=A_iX^2+B_iXY+C_iY^2$ having
these fixed coefficients $A_1$ and $A_2$. Recall (3.1). Then we have
$$N=O\left(2^{\omega\left(t^2-d_1d_2\right)}S\left(d_1,d_2,t^2\right
)\right)\eqno (3.17)$$
with implied absolute constant.

{\it Proof of Claim D.\/} Note that we can assume $A_i\neq 0$ for $
i=1,2$ by the assumption that $d_i$ is not a
square. Observe also that by the identity (3.6) with
$a=1,$ $b=0$ there are only two possibilities for $B_2A_1-B_1A_2$. So we can fix $
D$ such that there are at least $N/2$ inequivalent forms having the fixed coefficients
$A_1$, $A_2$ and having $B_2A_1-B_1A_2=D$. But we see by (3.15) that the residue of $
B_1$
modulo $2A_1$ determines the $SL_2({\bf Z})$-equivalence class of $\left
(Q_1,Q_2\right)$ once $A_1$, $A_2$ and $B_2A_1-B_1A_2$ is given. Therefore
it is enough to estimate the possible values of $B_1$ modulo
$A_1$ (for given $A_1$, $A_2$ and $D=B_2A_1-B_1A_2$) by the right-hand
side of (3.17).

Observe that we have
$$B_i^2\equiv d_i\pmod{A_i}\eqno (3.18)$$
for $i=1,2$. Let $p$ be a prime, and let us denote $\nu_p\left(A_
1\right)=\alpha$, $\nu_p\left(A_2\right)=\beta$. We consider two cases
separately.

(i) Assume $\nu_p\left(d_1\right)<\alpha$. Then we use (3.18) with
$i=1$. We see that there is a solution only if $\nu_p\left(d_1\right
)=2k$
for some integer $k$, and then we must have $\nu_p\left(B_1\right
)\ge k$
and $\left({{B_1}\over {p^k}}\right)^2\equiv{{d_1}\over {p^{2k}}}\hbox{\rm }\left
(p^{\alpha -2k}\right)$. Since $\alpha -2k>0$, we get from
this congruence that there are at most $2\left(1+\nu_p\left(2\right
)\right)$
possibilities for ${{B_1}\over {p^k}}$ modulo $p^{\alpha -2k}$. Hence we finally get
that there are at most $2\left(1+\nu_p\left(2\right)\right)p^{\left
[\nu_p\left(d_1\right)/2\right]}$
possibilities for $B_1$ modulo $p^{\alpha}$.

(ii) Assume $\nu_p\left(d_1\right)\ge\alpha$. Then we use again (3.18) with
$i=1$, and we see that we must have $\nu_p\left(B_1\right)\ge\alpha
/2$. So
there are at most $p^{\left[\alpha /2\right]}$ possibilities for $
B_1$ modulo $p^{\alpha}$.

Hence we see that in both cases there are at most $2\left(1+\nu_p\left
(2\right)\right)p^{\min\left(\left[{1\over 2}\nu_p\left(d_1\right
)\right],\left[{1\over 2}\alpha\right]\right)}$
possibilities for $B_1$ modulo $p^{\alpha}$. But we can see completely similarly that there are at most
$2\left(1+\nu_p\left(2\right)\right)p^{\min\left(\left[{1\over 2}
\nu_p\left(d_2\right)\right],\left[{1\over 2}\beta\right]\right)}$ possibilities for $
B_2$ modulo $p^{\beta}$.
Since $D=B_2A_1-B_1A_2$ is fixed, we see that if $B_2$ is given modulo $
p^{\beta}$,
then $D+B_1A_2=B_2A_1$ is given modulo $p^{\alpha +\beta}$, hence $
B_1$ is
given modulo $p^{\alpha}$. Taking into account (3.16) we finally get
that for every prime $p$ there are at most
$$2\left(1+\nu_p\left(2\right)\right)\min\left(p^{\left[\nu_p\left
(d_1\right)/2\right]},p^{\left[\nu_p\left(d_2\right)/2\right]},p^{\left
[\nu_p\left(t^2-d_1d_2\right)/2\right]}\right)$$
possibilities for $B_1$ modulo $p^{\nu_p\left(A_1\right)}$. We apply it for every
prime divisor $p$ of $A_1$ for which we have also that $p$
divides $t^2-d_1d_2$. If $p|A_1$ but $p$ does not divide $t^2-d_1
d_2$,
then by (3.16) we see that $p$ does not divide $A_2$, and so
$D=B_2A_1-B_1A_2$ implies that $B_1$ is determined modulo
$p^{\nu_p\left(A_1\right)}$. So for the number of possible values of $
B_1$
modulo $A_1$ we have the upper bound
$$C\left(\prod_{p|t^2-d_1d_2}2\right)\prod_p\min\left(p^{\left[\nu_
p\left(d_1\right)/2\right]},p^{\left[\nu_p\left(d_2\right)/2\right
]},p^{\left[\nu_p\left(t^2-d_1d_2\right)/2\right]}\right)\eqno (3
.19)$$
with an absolute constant $C$. Formula (3.19) proves Claim
D.

Claims C and D imply Lemma 3.1 at once.

{\bf 3.2. Upper bounds for certain special averages of }
$h\left(d_1,d_2,t\right)${\bf .} When we will apply Lemma 3.1, then we will have
numbers $d_i$ of special form $d_i=t^2_i-4$, and we will have
certain triple sums of $S\left(t^2_1-4,t^2_2-4,f^2\right)$, where $
t_1$, $t_2$, $f$
run over integers. We will use the trivial upper bounds
$S\left(t^2_1-4,t^2_2-4,f^2\right)\le S\left(t^2_1-4,t^2_2-4\right
)$ and
$S\left(t^2-4,t^2-4,f^2\right)\le S\left(t^2-4,f^2\right)$, and we will use
Lemmas 3.5, 3.6 and 3.7 below. We first need a preliminary lemma.

{\bf LEMMA 3.4.} {\it Let} $t_1,t_2>2$ {\it be integers,} $t_1\neq
t_2$, {\it and let} $E={\rm g}{\rm c}{\rm d}\left(t_1^2-4,t_2^2-4\right
)$.

{\it (i) There is a divisor} $e$ {\it of} $E$ {\it such that }
$$e\ge c\sqrt E,\quad e\left|t_1-\delta t_2\right.$$
{\it with an absolute constant} $c>0$ {\it and with some} $\delta
\in\left\{-1,1\right\}$.

{\it (ii) We have} $E\le\left|t_1-t_2\right|\left(t_1+t_2\right)$.

{\it Proof.\/} Part (ii) follows at once from the fact that $E$
divides $t_1^2-t_2^2$, so it remains to show part (i). Let $p|E$ be a
prime. Then $p^{\nu_p\left(E\right)}$ divides $\left(t_1-t_2\right
)\left(t_1+t_2\right)$, so writing $\alpha :=\nu_p\left(t_1-t_2\right
)$ and $\beta :=\nu_p\left(t_1+t_2\right)$ we have
$\nu_p\left(E\right)\le\alpha +\beta$. If $m=\min\left(\alpha ,\beta\right
)$, then $m\le\nu_p\left(2t_1\right)$, so
$2m\le\nu_p\left(4t_1^2\right)$. But $0<\nu_p\left(E\right)\le\nu_
p\left(4t_1^2-16\right)$, so if $m>0$, then
we must have $p=2$. If $p=2$ and $m>2$, then we have
$\nu_2\left(4t_1^2\right)>4$, and so $\nu_2\left(4t_1^2-16\right)
=4$, hence $\nu_2\left(E\right)\le 4$. It
follows for every prime $p$ that $p^{\nu_p\left(E\right)}$ divides either
$16\left(t_1-t_2\right)$ or $16\left(t_1+t_2\right)$. Then there is a decomposition
$E=e_1e_2$ such that ${\rm g}{\rm c}{\rm d}\left(e_1,e_2\right)=1$, and $
e_1$ divides
$16\left(t_1-t_2\right)$, $e_2$ divides $16\left(t_1+t_2\right)$. The lemma is proved.

{\bf LEMMA 3.5.} {\it Let} $3\le a<b\le c\le 2a$ {\it be integers. Recall the definition (3.1). For any} $
\epsilon >0$ {\it we have the following two inequalities:}
$$\sum_{t_1=a}^{c-1}\sum_{a\le t_2\le c-1,0<\left|t_2-t_1\right|\le
b-a}S\left(t_1^2-4,t_2^2-4\right)\ll_{\epsilon}a^{1/2+\epsilon}\left
(c-a\right)\left(b-a\right)^{1/2}\eqno (3.20)$$
{\it and}
$$\sum_{t_1=a}^{c-1}\sum_{a\le t_2\le c-1,0<\left|t_2-t_1\right|}{{
S\left(t_1^2-4,t_2^2-4\right)}\over {\sqrt {\left|t_1-t_2\right|}}}
\ll_{\epsilon}a^{1/2+\epsilon}\left(c-a\right).\eqno (3.21)$$
{\bf REMARK 3.2.} Estimating very crudely by
$S\left(t_1^2-4,t_2^2-4\right)\ll a$ the trivial bound in (3.20) would be
$\left(c-a\right)\left(b-a\right)a$, so in (3.20) we save roughly $\sqrt {\left
(b-a\right)a}$.

{\it Proof.\/} The second inequality follows at once from the
first one by usig a dyadic subdivision. To prove the
first one we may assume $b=c$, since the general case
follows from this special case by dividing the summation
over $t_1$ into $O\left({{c-a}\over {b-a}}\right)$ subsums.

So let $b=c$. Note that by Lemma 3.4 (ii) we have that
${\rm g}{\rm c}{\rm d}\left(t_1^2-4,t_2^2-4\right)\ll\left(b-a\right
)b.$ Then by Lemma 3.4 (i) and
by $S\left(t_1^2-4,t_2^2-4\right)\le\sqrt {{\rm g}{\rm c}{\rm d}\left
(t_1^2-4,t_2^2-4\right)}$ we have that the
left-hand side of (3.20) is
$$\ll\sum_{t_1=a}^{b-1}\left(\sum_{E|t_1^2-4,E\ll\left(b-a\right)
b}\sqrt E\sum_{e|E,e\ge c\sqrt E}\sum_{\delta\in\left\{-1,1\right
\}}\sum_{a\le t_2<b,e|t_2-\delta t_1}1\right),$$
and for a given $t_1$ the bracket is
$$\ll\sum_{E|t_1^2-4,E\ll\left(b-a\right)b}\sqrt E\sum_{e|E,e\ge
c\sqrt E}\left(1+{{b-a}\over e}\right)\ll b^{\epsilon}\left(\sqrt {\left
(b-a\right)b}+b-a\right).$$
The lemma follows.

{\bf LEMMA 3.6.} {\it Let} $3\le a<b\le 2a$ {\it be integers. Then for every} $
\epsilon >0$ {\it we have the following two inequalities\/}:
$$\sum_{t=a}^{b-1}\max\left\{k\ge 1:\;k^2|t^2-4\right\}\ll_{\epsilon}
a^{1+\epsilon}\sqrt {b-a}\eqno (3.22)$$
{\it and}
$$\sum_{t_1=a}^{b-1}\sum_{t_2=a}^{b-1}S\left(t_1^2-4,t_2^2-4\right
)\ll_{\epsilon}a^{\epsilon}\left(a\sqrt {b-a}+a^{1/2}\left(b-a\right
)^{3/2}\right).\eqno (3.23)$$
{\bf REMARK 3.3.} Estimating every summand by $O\left(a\right)$ in (3.22)
the trivial bound in (3.22) would be $\left(b-a\right)a$, so in (3.20)
we save roughly $\sqrt {b-a}$.

{\it Proof.\/} Statement (3.23) follows at once from Lemma 3.5
and (3.22), so we deal only with (3.22). It is enough to show that for any integer
$1\le K\le 2a$ we have that
$$K\sum_{t=a}^{b-1}\sum_{k=K}^{2K}\sum_{d>1,t^2-4=dk^2}\mu^2\left
(d\right)\ll_{\epsilon}a^{\epsilon}a\sqrt {b-a}.\eqno (3.24)$$
A trivial upper bound for the left-hand side of (3.24) is $K\left
(b-a\right)$.

Let $d$ be fixed and assume $t^2-4=dk^2.$ Then $\alpha :={{t+k\sqrt
d}\over 2}$ is an
algebraic integer, since it is a root of the equation
$x^2-tx+1$. We also see that $\alpha$ is a unit in the ring $R$ of
algebraic integers of the real quadratic field ${\bf Q}\left(\sqrt
d\right)$. By
the Dirichlet Unit Theorem there is a unit $1<\epsilon\in R$ such
that every unit of $R$ has the form $\pm\epsilon^l$ with integer $
l$.
One has $\epsilon ={{a+b\sqrt d}\over 2}$ with integers $a$, $b$, where $
b\neq 0$. Then
$\epsilon^{-1}=\delta{{a-b\sqrt d}\over 2}$ with $\delta\in\left\{
-1,1\right\}$, hence $\epsilon =b\sqrt d+\delta\epsilon^{-1}$, so
$\epsilon >\sqrt d-1\ge\sqrt 2-1$. But $\alpha =\epsilon^l$ with some positive integer $
l$ and
$\alpha\le t\le 2a$. So we proved that for a fixed $d$ there are at
most $C\log a$ possibilities for the pair $\left(t,k\right)$ with an
absolute constant $C$. We have $d\ll{{a^2}\over {K^2}}$ on the left-hand
side of (3.24), hence we finally get that the left-hand side of (3.24) is $
\ll_{\epsilon}{{a^{2+\epsilon}}\over K}$.

So the left-hand side of (3.24) is $\ll_{\epsilon}a^{\epsilon}\min\left
(K\left(b-a\right),{{a^2}\over K}\right)$. This minimum here is clearly $
\ll a\sqrt {b-a}$, and the lemma is proved.

{\bf LEMMA 3.7.} {\it Let} $t>2$ {\it be an integer and let} $1\le
A\ll t^2${\it . Then for any} $\epsilon >0$ {\it we have that}
$$\sum_{f\in {\bf Z},\,t^2-4-A\le\left|f\right|<t^2-4}{{S\left(t^
2-4,f^2\right)}\over {\sqrt {t^2-4-\left|f\right|}}}\ll_{\epsilon}
t^{\epsilon}\sqrt A.$$
{\bf REMARK 3.4.} We save roughly $t$ here with respect to the
trivial bound.

{\it Proof.\/} The left-hand side is at most
$$\sum_{k^2|t^2-4}k\sum_{g\in {\bf Z},\,0<{{t^2-4}\over k}-\left|
g\right|\le{A\over k}}{1\over {\sqrt k\sqrt {{{t^2-4}\over k}-\left
|g\right|}}},$$
and the inner sum here is $\ll{1\over {\sqrt k}}\sqrt {{A\over k}}
={{\sqrt A}\over k}$. We used here
that the inner sum is empty if $A<k$. The lemma is proved.

{\bf 4. Identities and estimates for special functions.}
\medskip

In this section we consider the functions ${\cal I}\left(t_1,t_2,
F,m_1,m_2\right)$ and ${\cal J}\left(t_1,t_2,m_1,m_2\right)$ defined in (2.10) and (2.9) for the special case when the functions $
m_i$ are characteristic functions of some intervals $\left[0,x_i\right
]$ for
$i=1,2$. This case will be important in our application. In
the first subsection we prove an identity for this ${\cal I}$-function for every $
1\neq F>0$, in the second and third subsections we use it to give estimates for the cases
$F>1$ and $F<1$, respectively. In the last subsection we
compute ${\cal J}\left(t_1,t_2,m_1,m_2\right)$ for the above-mentioned special case.

{\bf 4.1. Computing} ${\cal I}\left(t_1,t_2,F,m_1,m_2\right)$ {\bf when} $
m_i$ {\bf are characteristic functions.} For $S_0,T_0,F>0$, $F\neq
1$ define
$$Z\left(S_0,T_0,F\right):=\int\!\!\!\int{1\over {\sqrt {S^2+T^2+
2FTS+1-F^2}}}dSdT,\eqno (4.1)$$
where we integrate over the set
$$\left\{\left(S,T\right)\in {\bf R}^2:\;\left|S\right|\le S_0,\left
|T\right|\le T_0,\,S^2+T^2+2FTS+1-F^2>0\right\}.$$
By the reasoning of Remark 2.1 we see that (4.1) is absolutely convergent. One can see that (4.1) is divergent for $
F=1$,
but we do not need that case. It is clear that the ${\cal I}$-function can be expressed by the
$Z$-function in the case when $m_i$ are characteristic
functions.

{\bf LEMMA 4.1.} {\it Let} $S_0,T_0,F>0$, $F\neq 1$. {\it We have}
$$Z\left(S_0,T_0,F\right)=2J\left(S_0,T_0,F\right)+2J\left(T_0,S_
0,F\right),$$
{\it where we write}
$$J\left(S_0,T_0,F\right):=\int_{\left|y\right|\le T_0/S_0,\>\left
(1+y^2+2Fy\right)S_0^2>F^2-1}{{\sqrt {\left(1+y^2+2Fy\right)S_0^2
+1-F^2}}\over {1+y^2+2Fy}}dy\eqno (4.2)$$
{\it in the case} $F>1${\it , and }
$$J\left(S_0,T_0,F\right):=\int_{\left|y\right|\le T_0/S_0}{{\sqrt {\left
(1+y^2+2Fy\right)S_0^2+1-F^2}}\over {1+y^2+2Fy}}dy-{1\over 2}\int_{
-\infty}^{\infty}{{\sqrt {1-F^2}}\over {1+y^2+2Fy}}dy$$
{\it in the case} $F<1$.

{\it Proof.\/} It is clear by the substitution $\left(S,T\right)\rightarrow\left
(-S,-T\right)$ that the $S<0$ and $S>0$ parts of the
integral (4.1) have the same value. For $S>0$ we make the
substitution $y=T/S$, and we get
$$Z\left(S_0,T_0,F\right)=2\int_{-\infty}^{\infty}\int{S\over {\sqrt {
S^2+y^2S^2+2FyS^2+1-F^2}}}dSdy,$$
where the inner integral is taken over the set
$$\left\{S\in {\bf R}:\;0\le S\le S_0,\left|Sy\right|\le T_0,\,S^
2+y^2S^2+2FyS^2+1-F^2>0\right\}.\eqno (4.3)$$
If $F<1$, then the last condition is always true. If $F>1$,
then $1+y^2+2Fy>0$ should hold, otherwise (4.3) is empty.
For a fixed $y$ we integrate in $S$ over the interval
$$\sqrt {{{F^2-1}\over {1+y^2+2Fy}}}\le S\le\min\left(S_0,{{T_0}\over {\left
|y\right|}}\right)$$
in the case $F>1$, and we integrate over
$$0\le S\le\min\left(S_0,{{T_0}\over {\left|y\right|}}\right)$$
in the case $F<1$. We consider separately the cases
$\left|y\right|\le T_0/S_0$ and $\left|y\right|\ge T_0/S_0$. We can compute the
$S$-integral in every case. Making the substitution $y\rightarrow
1/y$
in the case $\left|y\right|\ge T_0/S_0$ we obtain the lemma.

{\bf 4.2. The case} $F>1${\bf .} In Lemma 4.3 we express the function
$J\left(S_0,T_0,F\right)$ defined in Lemma 4.1 in terms of a simple
function in the case $F>1$. Then we give estimates
for this simple function in Lemma 4.4 and for higher
derivatives of a variant of this function in Lemma 4.6.

{\bf LEMMA 4.2.} {\it Let} $S_0,T_0>0$ {\it and} $F>1$. {\it For the function }
$J\left(S_0,T_0,F\right)$ {\it defined in (4.2) we have}
$$J\left(S_0,T_0,F\right)=\int_{H\left(S_0,T_0,F\right)}{{\sqrt {\left
(1+y^2+2Fy\right)S_0^2+1-F^2}}\over {1+y^2+2Fy}}dy,\eqno (4.4)$$
{\it where the set} $H\left(S_0,T_0,F\right)$ {\it is defined as follows\/}:
$$H\left(S_0,T_0,F\right)=\left[-{{T_0}\over {S_0}},{{T_0}\over {
S_0}}\right],\hbox{\rm \ if }{{T_0}\over {S_0}}\le 1,\,1<F\le\,A\left
(S_0,T_0\right),\eqno (4.5)$$
$$H\left(S_0,T_0,F\right)=\left[-{{T_0}\over {S_0}},C\left(F,S_0\right
)\right]\cup\left[D\left(F,S_0\right),{{T_0}\over {S_0}}\right],\hbox{\rm \ if }{{
T_0}\over {S_0}}\ge 1,\,1<F\le\,A\left(S_0,T_0\right),\eqno (4.6)$$
$$H\left(S_0,T_0,F\right)=\left[D\left(F,S_0\right),{{T_0}\over {
S_0}}\right],\hbox{\rm \ if }\,A\left(S_0,T_0\right)\le F\le B\left
(S_0,T_0\right),\eqno (4.7)$$
$$H\left(S_0,T_0,F\right)=\emptyset ,\hbox{\rm \ if }\,F\ge B\left
(S_0,T_0\right),\eqno (4.8)$$
{\it where we write}
$$A\left(S_0,T_0\right):=\sqrt {\left(1+S_0^2\right)\left(1+T_0^2\right
)}-S_0T_0,\eqno (4.9)$$
$$B\left(S_0,T_0\right):=\sqrt {\left(1+S_0^2\right)\left(1+T_0^2\right
)}+S_0T_0,\eqno (4.10)$$
$$C\left(F,S_0\right):=-F-{{\sqrt {F^2-1}\sqrt {1+S_0^2}}\over {S_
0}},\;D\left(F,S_0\right):=-F+{{\sqrt {F^2-1}\sqrt {1+S_0^2}}\over {
S_0}}.$$
{\it We mean every statement in such a way that if we write an interval} $\left
[a,b\right]${\it , then this implicitly means that }
$a\le b$.

{\it Proof.\/} One can check that $\left(1+y^2+2Fy\right)S_0^2>F^
2-1$
holds if and only if $y<C\left(F,S_0\right)$ or $y>D\left(F,S_0\right
)$. The
following three claims can be checked by direct
computation. For the proof of Claim 2 we use the
obvious fact that $F-S_0T_0>-\sqrt {\left(1+S_0^2\right)\left(1+T_
0^2\right)}$.

{\bf CLAIM 1.} The sign of
$${{\sqrt {F^2-1}\sqrt {1+S_0^2}}\over {S_0}}-\left|{{T_0}\over {
S_0}}-F\right|$$
equals the sign of $F-A\left(S_0,T_0\right)$.

{\bf CLAIM 2.} The sign of
$${{\sqrt {F^2-1}\sqrt {1+S_0^2}}\over {S_0}}-\left({{T_0}\over {
S_0}}+F\right)$$
equals the sign of $F-B\left(S_0,T_0\right)$.

{\bf CLAIM 3.} The sign of ${{T_0}\over {S_0}}-1\,$ equals the sign of ${{
T_0}\over {S_0}}-A\left(S_0,T_0\right)$.

We get from Claim 1 that if ${{T_0}\over {S_0}}\le 1\,$ and $1<F\le\,
A\left(S_0,T_0\right)$,
then $D\left(F,S_0\right)\le -{{T_0}\over {S_0}}$, and this gives (4.5).

If ${{T_0}\over {S_0}}\ge 1\,$ and $1<F\le\,A\left(S_0,T_0\right)$, then we get from Claims 3 and 1 that
$$F+{{\sqrt {F^2-1}\sqrt {1+S_0^2}}\over {S_0}}\le{{T_0}\over {S_
0}},$$
and this implies (4.6).

If $A\left(S_0,T_0\right)\le F\le B\left(S_0,T_0\right)$, then $-{{
T_0}\over {S_0}}\le D\left(F,S_0\right)\le{{T_0}\over {S_0}}$ by
Claims 1 and 2, and $C\left(F,S_0\right)\le -{{T_0}\over {S_0}}$ by Claim 1. This proves (4.7).

If $F\ge B\left(S_0,T_0\right)$, then $D\left(F,S_0\right)\ge{{T_
0}\over {S_0}}$ by Claim 2, and
$C\left(F,S_0\right)\le -{{T_0}\over {S_0}}$ by Claim 1. This gives (4.8), and the lemma
is proved.

{\bf LEMMA 4.3.} {\it Use the notations of Lemma 4.2. Write} $\sigma
=1+{1\over {S_0^2}}${\it , and for} $0<y<1$ {\it let}
$$\Phi\left(y\right)=\Phi\left(S_0,y\right):=\int_0^y{{\sigma r^2}\over {\left
(1-r^2\right)\left(r^2+\sigma -1\right)}}dr.\eqno (4.11)$$
{\it Then for} $1<F\le B\left(S_0,T_0\right)$ {\it we have }
$J\left(S_0,T_0,F\right)=S_0\left(\Phi\left(y_1\right)+\epsilon\Phi\left
(y_2\right)\right),$ {\it where}
$$y_1=y_1\left(S_0,T_0,F\right):=\sqrt {1-{{\left(1+S_0^2\right)\left
(F^2-1\right)}\over {\left(T_0+S_0F\right)^2}}},$$
$$y_2=y_2\left(S_0,T_0,F\right):=\sqrt {1-{{\left(1+S_0^2\right)\left
(F^2-1\right)}\over {\left(T_0-S_0F\right)^2}}},$$
{\it and} $\epsilon =\epsilon\left(S_0,T_0,F\right)$ {\it is defined as follows:}
$$\epsilon =1,\hbox{\rm \ if }{{T_0}\over {S_0}}>1,\,1<F\le\,A\left
(S_0,T_0\right),$$
$$\epsilon =-1,\hbox{\rm \ if }{{T_0}\over {S_0}}<1,\,1<F\le\,A\left
(S_0,T_0\right),$$
$$\epsilon =0,\hbox{\rm \ if }\,A\left(S_0,T_0\right)<F\le B\left
(S_0,T_0\right).$$
{\it Assuming} $1<F\le B\left(S_0,T_0\right)$ {\it for} $i=1$ {\it and} $
1<F\le A\left(S_0,T_0\right)$ {\it for} $i=2$ {\it we have the following statement for} $
i=1,2$:
$$0\le y_i\le 1-c_1{{F-1}\over {\left(1+S_0\right)^{c_2}\left(1+T_
0\right)^{c_2}}}\eqno (4.12)$$
{\it with some positive absolute constants} $c_1$, $c_2$.

{\it Proof.\/} Note that $1<F\le\,A\left(S_0,T_0\right)$ implies ${{
T_0}\over {S_0}}\neq 1$ by
Claim 3, so $\epsilon$ is well-defined. We get the statement $0\le
y_i\le 1$ by Claims 1 and 2
above. Then (4.12) follows by easy estimates using the
conditions $S_0,T_0>0$ and $1<F\le B\left(S_0,T_0\right)$.

To compute $J\left(S_0,T_0,F\right)$ we use (4.4). Note that this is
the same formula as (4.2), but the integration set is given there explicitly. Use the substitution
$$r=r\left(y\right)=\sqrt {1-{{\sigma\left(F^2-1\right)}\over {\left
(y+F\right)^2}}}=\sqrt {{{1+y^2+2Fy+\left(F^2-1\right)\left(1-\sigma\right
)}\over {\left(y+F\right)^2}}}.$$
We have a positive number under the square root by
(4.2). It is clear by the conditions and by the definitions of
$C\left(F,S_0\right)$ and $D\left(F,S_0\right)$ that the sign of $
y+F$ is constant
on each of the four intervals which are present in (4.5),
(4.6) and (4.7), hence $r$ is well-defined and strictly monotone on each of these
intervals. It is easy to check that
$${{\sqrt {1+y^2+2Fy+\left(F^2-1\right)\left(1-\sigma\right)}}\over {
1+y^2+2Fy}}\left|{{dy}\over {dr}}\right|={{\sigma r^2}\over {\left
(1-r^2\right)\left(r^2+\sigma -1\right)}}.$$
We have $r\left(C\left(F,S_0\right)\right)=r\left(D\left(F,S_0\right
)\right)=0$, hence by Lemma
4.2 we get the present lemma.

{\bf LEMMA 4.4.} {\it Let} $S_0>0$ {\it and} $0<y<1${\it . Then for the function} $
\Phi\left(S_0,y\right)$ {\it defined in Lemma 4.3 we have the following estimates:}
$$S_0\Phi\left(S_0,y\right)\ll S_0^3y^3,\hbox{\rm \ if }S_0\ge 1,\;
0<y\le{1\over {2S_0}},\eqno (4.13)$$
$$S_0\Phi\left(S_0,y\right)\ll S_0y,\hbox{\rm \ if}\;S_0\ge 1,\;{
1\over {2S_0}}\le y\le{1\over 2},\eqno (4.14)$$
$$S_0\Phi\left(S_0,y\right)\ll S_0y^3,\hbox{\rm \ if}\;S_0\le 1,\;
y\le{1\over 2},\eqno (4.15)$$
{\it finally we have in every case that}
$$\Phi\left(S_0,y\right)\ll\log{1\over {1-y}}.\eqno (4.16)$$
{\it The implied constants are absolute in formulas (4.13)-(4.16).}

{\it Proof.\/} We have by the definitions that
$$S_0\Phi\left(S_0,y\right)=S_0\int_0^y{{\left(S_0^2r^2+r^2\right
)}\over {\left(1-r^2\right)\left(S_0^2r^2+1\right)}}dr.$$
Every estimate follows easily, the lemma is proved.

We recall Fa\`a di Bruno's formula. If $F$ and $G$ are smooth functions and $
H\left(x\right)=F\left(G\left(x\right)\right)$, then for every $j
\ge 1$ we have
$$H^{\left(j\right)}\left(x\right)=\sum_{l=1}^j\sum_{k=\left(k_1,
\ldots ,k_j\right)\in H_{j,l}}a_{j,l,k}F^{\left(l\right)}\left(G\left
(x\right)\right)\prod_{i=1}^j\left(G^{\left(i\right)}\left(x\right
)\right)^{k_i},\eqno (4.17)$$
with some constants $a_{j,l,k}$, where
$$H_{j,l}=\left\{\left(k_1,\ldots ,k_j\right)\in {\bf Z}^j:\;k_i\ge
0,\sum_{i=1}^jk_i=l,\sum_{i=1}^jik_i=j\right\}.\eqno (4.18)$$
This can be seen by induction using the chain rule.

{\bf LEMMA 4.5.} {\it Let} $S_0>0$, $\sigma =1+{1\over {S_0^2}}$.

{\it (i) Recall the definition of} $\Phi\left(y\right)$ {\it from (4.11). Write }
$\phi\left(t\right)=\Phi\left({1\over t}\right)${\it , then for every} $
j\ge 1$ {\it and} $t>1$ {\it we have }
$$\phi^{\left(j\right)}\left(t\right)\ll_j{1\over {t\left(t-1\right
)^j}}$$
{\it uniformly in} $S_0$.

{\it (ii) For} $Y>\sqrt {\sigma}$ {\it let} $G\left(Y\right)={Y\over {\sqrt {
Y^2-\sigma}}}${\it . Then for every} $j\ge 1$ {\it we have}
$$G^{\left(j\right)}\left(Y\right)\ll_j\left({Y\over {Y^2-\sigma}}\right
)^j{Y\over {\sqrt {Y^2-\sigma}}}\hbox{\rm \ for $\sqrt {\sigma}<Y
\le 2\sqrt {\sigma}$},\eqno (4.19)$$
$$G^{\left(j\right)}\left(Y\right)\ll_j{{\sigma}\over {Y^{j+2}}}\hbox{\rm \ for $
Y\ge 2\sqrt {\sigma}$},\eqno (4.20)$$
{\it uniformly in} $S_0$.

{\it (iii) For} $Y>\sqrt {\sigma}$ {\it let} $H\left(Y\right)=\phi\left
(G\left(Y\right)\right)$. {\it Then for every} $j\ge 1$ {\it and} $
Y>\sqrt {\sigma}$ {\it we have}
$$H^{\left(j\right)}\left(Y\right)\ll_j{{\sqrt {Y^2-\sigma}}\over
Y}\left({Y\over {Y^2-\sigma}}\right)^j$$
{\it uniformly in} $S_0$.

{\it Proof.\/} By (4.11) and the substitution $r\rightarrow 1/r$ we have
$$\phi\left(t\right)=\int_t^{\infty}{{\sigma}\over {\left(r^2-1\right
)\left(1+r^2\left(\sigma -1\right)\right)}}dr.$$
The integrand here equals ${1\over {r^2-1}}-$${1\over {r^2+S_0^2}}$. Considering
separately the cases $t\ge 2$ and $1<t\le 2$ we obtain (i) easily.

For the proof of (4.19) note that if $j\ge 1$, then $G^{\left(j\right
)}\left(Y\right)$ is
the linear combination of terms of the form ${{Y^l}\over {\left(\sqrt {
Y^2-\sigma}\right)^{j+l}}}$,
where $0\le l\le j+1$, and $j+l$ is odd. Here clearly $l=j+1$
gives the largest term, and we get (4.19). Statement (4.20)
follows easily from the Taylor expansion
$$G\left(Y\right)={1\over {\sqrt {1-\sigma Y^{-2}}}}=1+{{\sigma}\over {
2Y^2}}+\sum_{m=2}^{\infty}a_m{{\sigma^m}\over {Y^{2m}}},\eqno (4.
21)$$
where $a_m$ are absolute constants such that
$\sum_{m=1}^{\infty}\left|a_m\right|r^m<\infty$ for every $0<r<1$.

For the proof of (iii) we use Fa\`a di Bruno's formula (4.17), and we see that it is
enough to estimate terms of the form
$$\phi^{\left(l\right)}\left(G\left(Y\right)\right)\prod_{i=1}^j\left
(G^{\left(i\right)}\left(Y\right)\right)^{k_i},\eqno (4.22)$$
where $1\le l\le j$, and $k_i$ satisfy the conditions in (4.18).

If $Y\ge 2\sqrt {\sigma}$, then $1<G\left(Y\right)\le{2\over {\sqrt
3}}$, and we get from (i) and
(4.20) that (4.22) is
$$\ll_j{1\over {\left(G\left(Y\right)-1\right)^l}}\prod_{i=1}^j\left
({{\sigma}\over {Y^{i+2}}}\right)^{k_i}=\left({{\sigma}\over {Y^2\left
(G\left(Y\right)-1\right)}}\right)^lY^{-j},$$
where we used the conditions in (4.18). We see from
(4.21) that $1\ll{{\sigma}\over {Y^2\left(G\left(Y\right)-1\right
)}}\ll 1$, so we get (iii) for $Y\ge 2\sqrt {\sigma}$.

If $\sqrt {\sigma}<Y\le 2\sqrt {\sigma}$, then $G\left(Y\right)\ge{
2\over {\sqrt 3}}$, and we get from (i) and
(4.19) that (4.22) is
$$\ll_j{1\over {\left(G\left(Y\right)\right)^{l+1}}}\prod_{i=1}^j\left
(\left({Y\over {Y^2-\sigma}}\right)^i{Y\over {\sqrt {Y^2-\sigma}}}\right
)^{k_i}={1\over {G\left(Y\right)}}\left({{Y/G\left(Y\right)}\over {\sqrt {
Y^2-\sigma}}}\right)^l\left({Y\over {Y^2-\sigma}}\right)^j,$$
where we used the conditions in (4.18). By the definition
of $G\left(Y\right)$ we get (iii) also for this case. The lemma is
proved.

{\bf LEMMA 4.6.} {\it Let} $S_0>0${\it ,} $F>1${\it ,} $t\ge 3${\it ,} $
\tau\in\left\{-1,1\right\}$ {\it be given. If }
$x>t^2-4${\it , write  }
$$T_0=T_0\left(x\right):=\sqrt {{x\over {t^2-4}}-1}\eqno (4.23)$$
{\it and }
$$R\left(x\right):={{\left|\tau FS_0+T_0\left(x\right)\right|}\over {
S_0\sqrt {F^2-1}}}.\eqno (4.24)$$
{\it Let the number} $\sigma$ {\it and the function} $H$ {\it be defined as in Lemma 4.5, and let us define} $
K\left(x\right)=H\left(R\left(x\right)\right)$ {\it for} $x\in H_{
S_0,F,t,\tau}${\it , where }
$$\hbox{\rm $H_{S_0,F,t,1}:=\left\{x>t^2-4:\;F<\,B\left(S_0,T_0\left
(x\right)\right)\right\}$,}$$
$$\hbox{\rm $H_{S_0,F,t,-1}:=\left\{x>t^2-4:\;F<\,A\left(S_0,T_0\left
(x\right)\right)\right\}$}$$
{\it (see (4.9) and (4.10)). Then} $K$ {\it is well-defined}. {\it If} $
\tau =-1$,
{\it then} $K\left(x\right)$ {\it is a smooth function for} $\left
\{x\in H_{S_0,F,t,-1}:\;T_0\left(x\right)<S_0\right\}$ {\it and also for }
$\left\{x\in H_{S_0,F,t,-1}:\;T_0\left(x\right)>S_0\right\}$. {\it For every} $
j\ge 1$ {\it and every} $x$ {\it satisfying the above conditions we have }
$$K^{\left(j\right)}\left(x\right)\ll_j\left(x-t^2+4\right)^{-j}\max\left
(1,\left({{T_0\left|T_0-FS_0\right|}\over {\sqrt {\left(S_0^2+1\right
)\left(T_0^2+1\right)}\left(A\left(S_0,T_0\right)-F\right)}}\right
)^j\right)$$
{\it for} $\tau =-1${\it , and}
$$K^{\left(j\right)}\left(x\right)\ll_j\left(x-t^2+4\right)^{-j}\max\left
(1,\left({{T_0\left(T_0+FS_0\right)}\over {\sqrt {\left(S_0^2+1\right
)\left(T_0^2+1\right)}\left(B\left(S_0,T_0\right)-F\right)}}\right
)^j\right)$$
{\it for} $\tau =1$.

{\it Proof.\/} First note that we see by Claims 1,2 that
$R\left(x\right)>\sqrt {\sigma}$, hence $K\left(x\right)$ is well-defined. Note also that for
$\tau =-1$ we have $\left|-FS_0+T_0\right|=FS_0-T_0$ in the
case $S_0>T_0$, and $\left|-FS_0+T_0\right|=T_0-FS_0$ in the case
$S_0<T_0$. This follows by the conditions, using Claim 3. We
cannot have $T_0=S_0$ if $\tau =-1$, because $T_0=S_0$ implies
$A\left(S_0,T_0\right)=1$, so $1<F<A\left(S_0,T_0\right)$ is impossible.
Hence if  $\tau =-1$, then $R\left(x\right)$ is indeed a smooth function for
$T_0<S_0$, and also for $T_0>S_0$, so we can speak about the
derivatives of $K$.

We see from Fa\`a di Bruno's formula (4.17) that it is enough to estimate terms of the form
$$H^{\left(l\right)}\left(R\left(x\right)\right)\prod_{i=1}^j\left
(R^{\left(i\right)}\left(x\right)\right)^{k_i},\eqno (4.25)$$
where $1\le l\le j$, and $k_i$ satisfy the conditions in (4.18). It
is clear by the definitions that for $i\ge 1$ we have
$$R^{\left(i\right)}\left(x\right)\ll_i{{\left(x-t^2+4\right)^{{1\over
2}-i}}\over {S_0\sqrt {t^2-4}\sqrt {F^2-1}}}.$$
Hence, using also Lemma 4.5 (iii) we get that (4.25) is
$$\ll_j{{\sqrt {\left(R\left(x\right)\right)^2-\sigma}}\over {R\left
(x\right)}}\left({{R\left(x\right)}\over {\left(R\left(x\right)\right
)^2-\sigma}}\right)^l\prod_{i=1}^j\left({{\left(x-t^2+4\right)^{{
1\over 2}-i}}\over {S_0\sqrt {t^2-4}\sqrt {F^2-1}}}\right)^{k_i},$$
and by the conditions in (4.18), using also (4.23) this equals
$${{\sqrt {\left(R\left(x\right)\right)^2-\sigma}}\over {R\left(x\right
)}}\left({{R\left(x\right)T_0}\over {\left(\left(R\left(x\right)\right
)^2-\sigma\right)S_0\sqrt {F^2-1}}}\right)^l\left(x-t^2+4\right)^{
-j}.$$
Note that ${{\sqrt {\left(R\left(x\right)\right)^2-\sigma}}\over {
R\left(x\right)}}\le 1$, and it is easy to compute by
(4.24) and $\sigma =1+{1\over {S_0^2}}$ that
$$S_0\sqrt {F^2-1}{{\left(R\left(x\right)\right)^2-\sigma}\over {
R\left(x\right)}}={{\left(S_0^2+1\right)\left(T_0^2+1\right)-\left
(F-\tau S_0T_0\right)^2}\over {\left|\tau FS_0+T_0\right|}}.$$
The lemma is proved.

{\bf 4.3. The case} $F<1${\bf .} In Lemma 4.7 we give a new expression for the function
$J\left(S_0,T_0,F\right)$ defined in Lemma 4.1 in the case $F<1$, and we give also upper bounds for the new
expression. In Lemma 4.8 we give another new
expression for $J\left(S_0,T_0,F\right)$, expressing it in terms of a simple function. Then in Lemma 4.10 we give estimates for higher derivatives of a variant of this simple function.

{\bf LEMMA 4.7.} {\it Let} $0<S_0$, $0<T_0<T_0^{\ast}$, $0<F<1$. {\it We have}
$$J\left(S_0,T_0,F\right)+J\left(T_0,S_0,F\right)=K\left(S_0,T_0,
F\right)+K\left(T_0,S_0,F\right),\eqno (4.26)$$
{\it where}
$$K\left(S_0,T_0,F\right):=\int_{\left|y\right|\le T_0/S_0}{{\sqrt {\left
(1+y^2+2Fy\right)S_0^2+1-F^2}-\sqrt {1-F^2}}\over {1+y^2+2Fy}}dy.$$
{\it We also have that}
$$K\left(S_0,T_0,F\right)\ll{{S_0T_0}\over {\sqrt {1-F^2}}},\;K\left
(S_0,T_0^{\ast},F\right)-K\left(S_0,T_0,F\right)\ll{{S_0\left(T_0^{
\ast}-T_0\right)}\over {\sqrt {1-F^2}}}.\eqno (4.27)$$
{\it Proof.\/} The first statement follows from Lemma 4.1 and
$$\int_{\left|y\right|\le T_0/S_0}{{dy}\over {1+y^2+2Fy}}+\int_{\left
|y\right|\le S_0/T_0}{{dy}\over {1+y^2+2Fy}}=\int_{-\infty}^{\infty}{{
dy}\over {1+y^2+2Fy}},$$
which follows by the substitution $y\rightarrow 1/y$. We have that
$${{\sqrt {\left(1+y^2+2Fy\right)S_0^2+1-F^2}-\sqrt {1-F^2}}\over {
1+y^2+2Fy}}={{S_0^2}\over {\sqrt {\left(1+y^2+2Fy\right)S_0^2+1-F^
2}+\sqrt {1-F^2}}},$$
and this is $\le{{S_0^2}\over {\sqrt {1-F^2}}}$. The lemma follows.

{\bf LEMMA 4.8.} {\it Let} $S_0,T_0>0$ {\it and} $F<1$. {\it Write} $
\sigma =1+{1\over {S_0^2}}${\it , and for} $-1<t<1$ {\it let}
$$V\left(t\right)=V\left(S_0,t\right):=\int_0^t{{\sigma}\over {\left
(1-r^2\right)\left(1+\left(\sigma -1\right)r^2\right)}}dr.$$
{\it Then we have} $J\left(S_0,T_0,F\right)=S_0\left(V\left(s_1\right
)-V\left(s_2\right)\right)-{1\over 2}\int_{-\infty}^{\infty}{{\sqrt {
1-F^2}}\over {1+y^2+2Fy}}dy,$ {\it where}
$$s_1=s_1\left(S_0,T_0,F\right):={{S_0F+T_0}\over {\sqrt {\left(S_
0F+T_0\right)^2+\left(1+S_0^2\right)\left(1-F^2\right)}}},$$
$$s_2=s_2\left(S_0,T_0,F\right):={{S_0F-T_0}\over {\sqrt {\left(S_
0F-T_0\right)^2+\left(1+S_0^2\right)\left(1-F^2\right)}}}.$$
{\it Proof}. We use the substitution
$$r=r\left(y\right)={{y+F}\over {\sqrt {\left(y+F\right)^2+\sigma\left
(1-F^2\right)}}}$$
in the first integral in the definition of $J\left(S_0,T_0,F\right
)$. It is easy to check that
$${{\sqrt {1+y^2+2Fy+\left(1-F^2\right)\left(\sigma -1\right)}}\over {
1+y^2+2Fy}}dy={{\sigma}\over {\left(1-r^2\right)\left(1+\left(\sigma
-1\right)r^2\right)}}dr,$$
and the lemma follows.

{\bf LEMMA 4.9.} {\it Let} $S_0>0$, $\sigma =1+{1\over {S_0^2}}$. {\it Recall the function} $
V$ {\it from Lemma 4.8.}

{\it (i) For every} $j\ge 1$ {\it and} $-1<t<1$ {\it we have that }
$$V^{\left(j\right)}\left(t\right)\ll_j\left(1-\left|t\right|\right
)\left({1\over {\left|t\right|+S_0}}+{1\over {1-\left|t\right|}}\right
)^{j+1}$$
{\it uniformly in} $S_0$.

{\it (ii) For} $-\infty <Y<\infty$ {\it let} $g\left(Y\right)={Y\over {\sqrt {
Y^2+\sigma}}}${\it . Then for every} $j\ge 1$ {\it we have}
$$g^{\left(j\right)}\left(Y\right)\ll_j\left({1\over {\sqrt {\sigma}}}\right
)^j\hbox{\rm \ for $\left|Y\right|\le 2\sqrt {\sigma}$},\eqno (4.
28)$$
$$g^{\left(j\right)}\left(Y\right)\ll_j{{\sigma}\over {\left|Y\right
|^{j+2}}}\hbox{\rm \ for $\left|Y\right|\ge 2\sqrt {\sigma}$},\eqno
(4.29)$$
{\it uniformly in} $S_0$.

{\it (iii) For} $-\infty <Y<\infty$ {\it let} $h\left(Y\right)=V\left
(g\left(Y\right)\right)$. T{\it hen for every} $j\ge 1$
{\it and} $Y>0$ {\it we have}
$$h^{\left(j\right)}\left(Y\right)\ll_j\sqrt {\sigma}\left(1+\left
|Y\right|\right)^{-j}$$
{\it uniformly in} $S_0$.

{\it Proof.\/} For the proof of (i) note that
$${{\sigma}\over {\left(1-r^2\right)\left(1+\left(\sigma -1\right
)r^2\right)}}={1\over {1-r^2}}+{1\over {S_0^2+r^2}}={1\over {1-r^
2}}+{{1/2S_0}\over {S_0+ir}}+{{1/2S_0}\over {S_0-ir}}.$$
Considering first the case $\left|t\right|\ge{1\over 2}$, and then if $\left
|t\right|\le{1\over 2}$, then considering separately $\left|t\right
|\le S_0$ and $\left|t\right|\ge S_0$ we get (i). In (ii) we can assume $
Y\ge 0$, and then the proof is completely similar to the
proof of (ii) of Lemma 4.5.

For the proof of (iii) we use Fa\`a di Bruno's formula (4.17), and we see that it is
enough to estimate terms of the form
$$V^{\left(l\right)}\left(g\left(Y\right)\right)\prod_{i=1}^j\left
(g^{\left(i\right)}\left(Y\right)\right)^{k_i},\eqno (4.30)$$
where $1\le l\le j$, and $k_i$ satisfy the conditions in (4.18).

If $\left|Y\right|\ge 2\sqrt {\sigma}$, then $\left|g\left(Y\right
)\right|\ge{2\over {\sqrt 5}}$, and we get from (i) and (4.29) that (4.30) is
$$\ll_j{1\over {\left(1-\left|g\left(Y\right)\right|\right)^l}}\prod_{
i=1}^j\left({{\sigma}\over {\left|Y\right|^{i+2}}}\right)^{k_i}=\left
({{\sigma}\over {\left|Y\right|^2\left(1-\left|g\left(Y\right)\right
|\right)}}\right)^l\left|Y\right|^{-j},$$
where we used the conditions in (4.18). It is easy to see
that $1\ll{{\sigma}\over {\left|Y\right|^2\left(1-\left|g\left(Y\right
)\right|\right)}}\ll 1$, so taking into account $\sigma\ge 1$
we get (iii) for $\left|Y\right|\ge 2\sqrt {\sigma}$. If $\left|Y\right
|\le 2\sqrt {\sigma}$, then $\left|g\left(Y\right)\right|\le{2\over {\sqrt
5}}$, and we get from (i) and
(4.28) that (4.30) is
$$\ll_j\left({1\over {\left|g\left(Y\right)\right|+S_0}}+1\right)^{
l+1}\prod_{i=1}^j\left(\left({1\over {\sqrt {\sigma}}}\right)^i\right
)^{k_i}=\left({1\over {\left|g\left(Y\right)\right|+S_0}}+1\right
)^{l+1}\left({1\over {\sqrt {\sigma}}}\right)^j,$$
where we used the conditions in (4.18). If $S_0\gg 1$, then $1\ll
\sigma\ll 1$, and we get (iii). If $S_0\ll 1$, then ${1\over {\sqrt {
\sigma}}}\ll S_0\ll{1\over {\sqrt {\sigma}}}$,
${{\left|Y\right|}\over {\sqrt {\sigma}}}\ll\left|g\left(Y\right)\right
|\ll{{\left|Y\right|}\over {\sqrt {\sigma}}}$, and
$$\left({1\over {\left|g\left(Y\right)\right|+S_0}}+1\right)^{l+1}\left
({1\over {\sqrt {\sigma}}}\right)^j\ll\left({1\over {\left|g\left
(Y\right)\right|+S_0}}\right)\left({1\over {\left|g\left(Y\right)\right
|\sqrt {\sigma}+S_0\sqrt {\sigma}}}\right)^j.$$
The lemma follows.

{\bf LEMMA 4.10.} {\it Let} $S_0>0${\it ,} $F<1${\it ,} $t\ge 3${\it ,} $
\tau\in\left\{-1,1\right\}$ {\it be given. If }
$x>t^2-4${\it , write  }
$$T_0=T_0\left(x\right):=\sqrt {{x\over {t^2-4}}-1}\eqno (4.31)$$
{\it and }
$$r\left(x\right):={{FS_0+\tau T_0}\over {S_0\sqrt {1-F^2}}}.\eqno
(4.32)$$
{\it Let the number} $\sigma$ {\it and the function} $h$ {\it be defined as in Lemma 4.9, and let us define} $
k\left(x\right)=h\left(r\left(x\right)\right)$ {\it for every} $x$ {\it satisfying} $
x>t^2-4$. {\it Then for every} $j\ge 1$ {\it and every} $x>t^2-4$ {\it we have}
$$k^{\left(j\right)}\left(x\right)\ll_j\sqrt \sigma\left(x-t^2+4\right
)^{-j}\max\left(1,\left({{T_0}\over {\left(S_0\sqrt {1-F^2}+\left
|\tau FS_0+T_0\right|\right)}}\right)^j\right).$$
{\it Proof.\/} We see from Fa\`a di Bruno's formula (4.17) that it is enough to estimate terms of the form
$$h^{\left(l\right)}\left(r\left(x\right)\right)\prod_{i=1}^j\left
(r^{\left(i\right)}\left(x\right)\right)^{k_i},\eqno (4.33)$$
where $1\le l\le j$, and $k_i$ satisfy the conditions in (4.18). It is clear by the definitions that for $
i\ge 1$ we have
$$r^{\left(i\right)}\left(x\right)\ll_i{{\left(x-t^2+4\right)^{{1\over
2}-i}}\over {S_0\sqrt {t^2-4}\sqrt {1-F^2}}}.$$
Hence, using also Lemma 4.9 (iii) we get that (4.33) is
$$\ll_j\sqrt {\sigma}\left(1+\left|r\left(x\right)\right|\right)^{
-l}\prod_{i=1}^j\left({{\left(x-t^2+4\right)^{{1\over 2}-i}}\over {
S_0\sqrt {t^2-4}\sqrt {1-F^2}}}\right)^{k_i}.$$
Using the conditions in (4.18) and the relations (4.31), (4.32) we obtain the lemma.

{\bf 4.4. Computing} ${\cal J}\left(t_1,t_2,m_1,m_2\right)$ {\bf when} $
m_i$ {\bf are characteristic functions.} For $x>0$ introduce the notation
$$k_x\left(y\right)=1\hbox{\rm \ for $0\le y\le x$},\quad k_x\left
(y\right)=0\hbox{\rm \ for $y>x$}.\eqno (4.34)$$
{\bf LEMMA 4.11.} {\it Let} $t_i>2$ {\it and} $x_i>0$ {\it for} $
i=1,2$. {\it Then} ${\cal J}\left(t_1,t_2,k_{x_1/4},k_{x_2/4}\right)$ {\it is nonzero only if} $
x_i>t_i^2-4$ {\it for} $i=1,2$. {\it Assuming that this is true, we have}
$${\cal J}\left(t_1,t_2,k_{x_1/4},k_{x_2/4}\right)=2{{\sqrt {1-m}}\over {\sqrt
m}},$$
{\it where} $m:=\max\left({{t_1^2-4}\over {x_1}},{{t_2^2-4}\over {
x_2}}\right)$.

{\it Proof.\/} The statement is trivial for $m\ge 1$, so let us
assume $m<1$. Then by definition we
have
$${\cal J}\left(t_1,t_2,k_{x_1/4},k_{x_2/4}\right)=\int_{-\arccos\left
(\sqrt m\right)}^{\arccos\left(\sqrt m\right)}{{d\theta}\over {\cos^
2\theta}}=2{{\sin\left(\arccos\left(\sqrt m\right)\right)}\over {\cos\left
(\arccos\left(\sqrt m\right)\right)}}.$$
The lemma is proved.

{\bf 5. First steps of the proof of Theorem 1.1. }
\medskip

We introduce some notation. Let $\eta_0$ be a given nonnegative smooth function on $\left
(0,\infty\right)$ such that $\eta_0\left(\tau\right)=0$
for $\tau\notin\left[1,2\right]$, and
$$\int_1^2\eta_0\left(\tau\right)d\tau =1.\eqno (5.1)$$
Recall the definition of $k_x\left(y\right)$ from (4.34). For $x>
0$, $D>0$, define
$$k_{x,D}\left(y\right):={1\over D}\int_D^{2D}\eta_0\left({{\tau}\over
D}\right)k_x\left(y+\tau\right)d\tau\eqno (5.2)$$
for $y\ge 0$. We will use also the notations of Theorem 1.1.

{\bf 5.1. A spectral estimate.} Our aim is to prove Lemma 5.2,
whose result will show that for a smoothed version of the
hyperbolic circle problem one can give a good estimate by
spectral methods. We first need some notation.

The hyperbolic Laplace operator is denoted by $\Delta :=y^2\left({{
\partial^2}\over {\partial x^2}}+{{\partial^2}\over {\partial y^2}}\right
)$. Let $\left\{u_j(z):\;j\ge 0\right\}$ be a complete orthonormal system of Maass forms for $
PSL_2({\bf Z})$ (the function $u_0(z)$ is constant), let
$\Delta u_j$$=\left(-{1\over 4}-t_j^2\right)u_j$, where $t_0={i\over
2}$ and $t_j$ is real for $j>0$.

If $m$ is a compactly supported bounded function on $[0,\infty )$, let (see [I], (1.62))
$$g_m\left(a\right)=2q_m\left({{e^a+e^{-a}-2}\over 4}\right),\hbox{\rm \ where $
q_m$$\left(v\right)=\int_0^{\infty}{{m\left(v+\tau\right)}\over {\sqrt {
\tau}}}d\tau$},\eqno (5.3)$$
and for any complex $r$ let
$$\hbox{\rm $h_m\left(r\right)=\int_{-\infty}^{\infty}g_m\left(a\right
)e^{ira}da.$}\eqno (5.4)$$
For simplicity introduce the abbreviations $h_x=h_{k_x}$ and
$h_{x,D}=h_{k_{x,D}}$ (see (5.3) and (5.4)).

{\bf LEMMA 5.1.} {\it Assume} $1<D<x/10$. {\it For every integer} $
j\ge 0$ {\it we have for} $r\ge 1$ {\it that}
$$h_{x,D}\left(r\right)\ll_j{{x^{1/2}}\over {r^{3/2}}}\min\left
(1,{x\over {Dr}}\right)^j+{{x^{1/2}}\over {r^{5/2}}}.\eqno (5.5)$$
{\it We also have for every real} $r$ {\it that}
$$h_{x,D}\left(r\right)\ll x^{{1\over 2}}\log x.\eqno
(5.6)$$
{\it Furthermore, we have}
$$h_{x,D}\left({i\over 2}\right)=4\pi x-4\pi D\int_1^2\eta_0\left
(\tau\right)\tau d\tau .\eqno (5.7)$$
{\it Proof.\/} It is easy to see by (5.3), (5.4) and (5.2) that
$$h_{x,D}\left(r\right)={1\over D}\int_D^{2D}\eta_0\left({{\tau}\over
D}\right)h_{x-\tau}\left(r\right)d\tau\eqno (5.8)$$
for every complex $r$. Now we apply Lemma 2.4 of [C]
for the function $h_{x-\tau}\left(r\right)$ choosing $R=R$$\left(
\tau\right)$ in that lemma in
such a way that
$${{\cosh\hbox{\rm $R$$\left(\tau\right)$}}\over 2}-{1\over 2}=x-
\tau\eqno (5.9)$$
holds. Applying part (d) of that lemma we
see that $h_{x-\tau}\left({i\over 2}\right)=4\pi\left(x-\tau\right
)$, and taking into account
(5.1) we get (5.7). The estimate (5.6) follows from a trivial
estimation of (2.6) of [C]. Finally, for the proof of (5.5) we apply part (a) of
Lemma 2.4 of [C]. Applying it in (5.8) we get for $r\ge 1$
that
$$h_{x,D}\left(r\right)={{2\sqrt {2\pi}}\over {r^{3/2}D}}\int_D^{
2D}\eta_0\left({{\tau}\over D}\right)\sqrt {\sinh\hbox{\rm $R$$\left
(\tau\right)$}}\cos\left(r\hbox{\rm $R$$\left(\tau\right)$}-{{3\pi}\over
4}\right)d\tau +O\left({{x^{1/2}}\over {r^{5/2}}}\right).$$
By the substitution $R=R\left(\tau\right)$, using (5.9) this equals
$${{\sqrt {2\pi}}\over {r^{3/2}D}}\int_{R_1}^{R_2}\eta_0\left({{1
+2x-\cosh\hbox{\rm $R$}}\over {2D}}\right)\left(\sinh\hbox{\rm $R$}\right
)^{3/2}\cos\left(r\hbox{\rm $R$}-{{3\pi}\over 4}\right)dR+O\left({{
x^{1/2}}\over {r^{5/2}}}\right),$$
where $\cosh\hbox{\rm $R$$_1$}=1+2x-4D,\;\cosh\hbox{\rm $R$$_2$}=
1+2x-2D.$ Repeated partial integration gives (5.5). The lemma is proved.

{\bf LEMMA 5.2.} {\it Assume} $1<D<x/10$ {\it and} $z\in\Omega${\it . Then}
$$\sum_{\gamma\in\Gamma}k_{x,D}\left(u\left(\gamma z,z\right)\right
)=12x-12D\int_1^2\eta_0\left(\tau\right)\tau d\tau +O_{\Omega}\left
({x\over {\sqrt D}}+x^{1/2}\log x\right).\eqno (5.10)$$
{\it Proof.\/} It is clear by (5.3) and (5.4) that the function $
h_{x,D}\left(r\right)$ satisfies condition (1.63) of [I], i.e it is even, it is holomorphic in the
strip $\left|\hbox{\rm Im}r\right|\le{1\over 2}+$$\epsilon$ and $
h_{k_{x,D}}\left(r\right)=O\left(\left(1+\left|r\right|\right)^{-
2-\epsilon}\right)$ in this
strip for some $\epsilon >0${\it .\/} Then we get from Theorem 7.4 of
[I], using again the abbreviation $h_{x,D}=h_{k_{x,D}}$ that the
left-hand side of (5.10) equals
$$\sum_{j=0}^{\infty}h_{x,D}\left(t_j\right)\left|u_j\left(z\right
)\right|^2+{1\over {4\pi}}\int_{-\infty}^{\infty}h_{x,D}\left(r\right
)\left|E\left(z,{1\over 2}+ir\right)\right|^2dr$$
for any $z\in \bbb H$, where $E\left(z,s\right)$ is the Eisenstein series for
$\Gamma =PSL_2({\bf Z})$, see [I], Chapter 3. Since for the
fundamental domain ${\cal F}$ defined in (1.3) we have that $\left
|{\cal F}\right|={{\pi}\over 3}$
(see [I], (6.33) and (3.26)), $\left|u_0\left(z\right)\right|^2$ equals ${
3\over {\pi}}$ for every
$z$. Then by Lemma 5.1 above and by [I], Proposition 7.2 we get the
lemma.

{\bf 5.2. Nonhyperbolic elements.} We give an easy estimate for the contribution of the nonhyperbolic elements in the hyperbolic circle problems.

{\bf LEMMA 5.3.} {\it Let} $z\in\Omega$ {\it and} $X>2${\it . Then for every} $
\epsilon >0$ {\it we have that}
$$\left|\left\{\gamma\in PSL_2({\bf Z}):\hbox{\rm \ $\left|\hbox{\rm $\hbox{\rm tr}\left
(\gamma\right)$}\right|\le 2,\;$}4u\left(\gamma z,z\right)\le X-2\right
\}\right|\ll_{\Omega ,\epsilon}X^{{1\over 2}+\epsilon}.\eqno (5.1
1)$$
{\it Proof.\/} Write $\gamma =\left(\matrix{a&b\cr
c&d\cr}
\right)$. First note that by [I], (1.9) and
(1.11) we have
$$4u\left(\gamma z,z\right)={{\left|cz^2+\left(d-a\right)z-b\right
|^2}\over {\hbox{\rm Im$^2z$}}}.$$
It is easy to compute that if $z=x+iy$, then we have
$$\hbox{\rm Im}\left(cz^2+\left(d-a\right)z-b\right)=2cxy+\left(d
-a\right)y,$$
$$\hbox{\rm Re}\left(cz^2+\left(d-a\right)z-b\right)=c\left(x^2-y^
2\right)+\left(d-a\right)x-b.$$
Hence if $z\in\Omega$, then the second inequality in (5.11) gives
$$2cx+d-a\ll_{\Omega}\sqrt X,\eqno (5.12)$$
$$c\left(x^2+y^2\right)+b\ll_{\Omega}\sqrt X.\eqno (5.13)$$
By the first inequality in (5.11) we get from (5.12) that
$$d=-cx+O_{\Omega}\left(\sqrt X\right),\;a=cx+O_{\Omega}\left(\sqrt
X\right),$$
and from these relations and (5.13) we get
$$1=ad-bc=-c^2x^2+c^2\left(x^2+y^2\right)+O_{\Omega}\left(\sqrt X\left
(\sqrt X+\left|c\right|\right)\right).$$
This implies $c=O_{\Omega}\left(\sqrt X\right)$, and so (5.12) gives
$d-a\ll_{\Omega}\sqrt X.$ Then there are $O_{\Omega}\left(\sqrt X\right
)$ possibilities for the pair $\left(a,d\right)$.
If $a$ and $d$ are given with $ad\neq 1$, then $bc=ad-1$ implies
that there are $O_{\epsilon}\left(X^{\epsilon}\right)$ possibilities for the pair $\left
(b,c\right)$.
Finally, if $ad=1$, then $bc=0$, and so (5.13) implies that
there are $\ll_{\Omega}\sqrt X$ possibilities for the pair $\left
(b,c\right)$. The
lemma is proved.

{\bf 5.3. Reduction to the estimation of a square integral on the fundamental domain.} For simplicity let us write
$N\left(z,X\right)=N\left(z,z,X\right)$. Let us take an integer $
J\ge 2$, it will
be fixed but we will choose it sufficiently large. Let $d$
be a parameter that will be chosen optimally later, at
the moment we assume that $d\ge 100$ and $100Jd\le X$.

Let us define
$$N_{d,J}\left(z,X\right):=\sum_{j=0}^J\left(-1\right)^j\left(\matrix{
J\cr
j\cr}
\right)\int_1^2\eta_0\left(\tau\right)N\left(z,X-jd\tau\right)d\tau
.\eqno (5.14)$$
Then using (5.1) we see that $N_{d,J}\left(z,X\right)$ equals
$$N\left(z,X\right)+\sum_{j=1}^J\left(-1\right)^j\left(\matrix{J\cr
j\cr}
\right)\sum_{\gamma\in\Gamma}\int_1^2\eta_0\left(\tau\right)k_{\left
(X-2\right)/4}\left(u\left(\gamma z,z\right)+{{jd\tau}\over 4}\right
)d\tau ,$$
which equals
$$N\left(z,X\right)+\sum_{j=1}^J\left(-1\right)^j\left(\matrix{J\cr
j\cr}
\right)\sum_{\gamma\in\Gamma}k_{\left(X-2\right)/4,jd/4}\left(u\left
(\gamma z,z\right)\right)$$
by (5.2). Applying Lemma 5.2 this equals
$$N\left(z,X\right)+O_{\Omega ,J}\left({X\over {\sqrt d}}+X^{1/2}\log
X\right)+\sum_{j=1}^J\left(-1\right)^j\left(\matrix{J\cr
j\cr}
\right)\left(3X-3jd\int_1^2\eta_0\left(\tau\right)\tau d\tau\right
)$$
for $z\in\Omega$. Now, $\sum_{j=1}^J\left(-1\right)^j\left(\matrix{
J\cr
j\cr}
\right)=-1,\;\sum_{j=1}^J\left(-1\right)^jj\left(\matrix{J\cr
j\cr}
\right)=0,$ which follows from the binomial theorem taking into
account that $j\left(\matrix{J\cr
j\cr}
\right)=J$$\left(\matrix{J-1\cr
j-1\cr}
\right)$ for $1\le j\le J$. Hence we
proved that for $z\in\Omega$ we have
$$N_{d,J}\left(z,X\right)=N\left(z,X\right)-3X+O_{\Omega ,J}\left
({X\over {\sqrt d}}+X^{1/2}\log X\right).\eqno (5.15)$$
Recalling the notation $M_{t,m}$ from (2.8) we get from Lemma 5.3 that
$$N\left(z,X\right)=O_{\Omega ,\epsilon}\left(X^{{1\over 2}+\epsilon}\right
)+\sum_{t>2}M_{t,k_{\left(X-2\right)/4}}\left(z\right)$$
for $z\in\Omega$, $X>2$ and for any $\epsilon >0$. Hence by (5.14) we see that
$$N_{d,J}\left(z,X\right)=O_{\Omega ,\epsilon ,J}\left(X^{{1\over
2}+\epsilon}\right)+\int_1^2\eta_0\left(\tau\right)\left(\sum_{t>
2}\sum_{j=0}^J\left(-1\right)^j\left(\matrix{J\cr
j\cr}
\right)M_{t,k_{\left(X-jd\tau -2\right)/4}}\left(z\right)\right)d
\tau$$
for $z\in\Omega$, $\epsilon >0$. By Cauchy-Schwarz we have that
$$\left(\int_1^2\eta_0\left(\tau\right)\left(\sum_{t>2}\sum_{j=0}^
J\left(-1\right)^j\left(\matrix{J\cr
j\cr}
\right)M_{t,k_{\left(X-jd\tau -2\right)/4}}\left(z\right)\right)d
\tau\right)^2$$
is
$$\ll\int_1^2\left(\sum_{t>2}\sum_{j=0}^J\left(-1\right)^j\left(\matrix{
J\cr
j\cr}
\right)M_{t,k_{\left(X-jd\tau -2\right)/4}}\left(z\right)\right)^
2d\tau .$$
Hence, using also (5.15) we finally get that
$$\int_{\Omega}\left(N\left(z,X\right)-3X\right)^2d\mu_z\eqno (5.
16)$$
is
$$O_{\Omega ,\epsilon ,J}\left({{X^2}\over d}+X^{1+\epsilon}+\int_
1^2\int_{{\cal F}}\left(\sum_{t>2}\sum_{j=0}^J\left(-1\right)^j\left
(\matrix{J\cr
j\cr}
\right)M_{t,k_{\left(X-jd\tau -2\right)/4}}\left(z\right)\right)^
2d\mu_zd\tau\right)\eqno (5.17)$$
if $\epsilon >0$, $d\ge 100$ and $100Jd\le X$.

We will show in Section 6 that if $\epsilon >0$ is given and the
integer $J\ge 2$ is fixed to be large enough in terms of $\epsilon$,
then we have
$$_{}\int_{{\cal F}}\left(\sum_{t>2}\sum_{j=0}^J\left(-1\right)^j\left
(\matrix{J\cr
j\cr}
\right)M_{t,k_{\left(X-jd\tau -2\right)/4}}\left(z\right)\right)^
2d\mu_z\ll_{\epsilon}X^{\epsilon}{{d^{5/2}}\over {\sqrt X}}\eqno
(5.18)$$
uniformly for $1\le\tau\le 2$ and
$$X^{2/3}\le d\le X^{99/100}.\eqno (5.19)$$
Assume that (5.18) is true. Then we see from (5.16) and
(5.17) that (5.16) equals
$$O_{\Omega ,\epsilon}\left({{X^2}\over d}+X^{\epsilon}{{d^{5/2}}\over {\sqrt
X}}\right)\eqno (5.20)$$
for any $d$ satisfying (5.19). Note that we choose $J$ in terms of $
\epsilon$, so we do not have
to denote the dependence on $J$ in (5.20). Choosing $d=X^{5/7}$
we obtain Theorem 1.1. So it is enough to show the estimate (5.18).

{\bf 6. Conclusion. }
\medskip

The goal of this section is to prove the estimate (5.18).

{\bf 6.1. Application of Lemma 2.2 and basic observations.} It is easy to see that if $
\gamma\in SL_2({\bf R})$ and the trace of $\gamma$
is $t>2$, then we have $u\left(\gamma z,z\right)\ge{{t^2-4}\over
4}$ for every $z\in\bbb H$. Therefore, the contribution of the terms
$t>\sqrt {X+2}$ to the sum (5.18) is $0$. We can take integers $1
\le I\ll\log X$ and
$$3=a_1<a_2<\ldots <a_I<1+\sqrt {X+2}\le a_{I+1}<2+\sqrt {X+2}\eqno
(6.1)$$
such that
$$a_{i+1}\le{3\over 2}a_i,\;3+\sqrt {X+2}-a_i\le 2\left(3+\sqrt {
X+2}-a_{i+1}\right)\eqno (6.2)$$
for $1\le i\le I$. By the Cauchy-Schwarz inequality we get for every $
1\le\tau\le 2$ that
$$\int_{\cal F}\left(\sum_{t>2}\sum_{j=0}^J\left(-1\right)^j\left(\matrix{
J\cr
j\cr}
\right)M_{t,k_{\left(X-jd\tau -2\right)/4}}\left(z\right)\right)^
2d\mu_z\ll\log X\sum_{i=1}^IU_i\eqno (6.3)$$
with
$$U_i=U_i\left(\tau\right):=\int_{\cal F}\left(\sum_{t=a_i}^{a_{i+1}-1}\sum_{
j=0}^J\left(-1\right)^j\left(\matrix{J\cr
j\cr}
\right)M_{t,k_{\left(X-jd\tau -2\right)/4}}\left(z\right)\right)^
2d\mu_z.\eqno (6.4)$$
By Lemma 2.2 we have for every $1\le i\le I$ that $U_i$ equals
the sum of
$$\sum_{t_1=a_i}^{a_{i+1}-1}\sum_{t_2=a_i}^{a_{i+1}-1}E_{t_1,t_2}
S_{t_1,t_2},\eqno (6.5)$$
$$\sum_{t_1=a_i}^{a_{i+1}-1}\sum_{t_2=a_i}^{a_{i+1}-1}\sum_{f\in
{\bf Z},f^2<\left(t_1^2-4\right)\left(t_2^2-4\right)}h\left(t_1^2
-4,t_2^2-4,f\right)R_{t_1,t_2,f}\eqno (6.6)$$
and
$$\sum_{t_1=a_i}^{a_{i+1}-1}\sum_{t_2=a_i}^{a_{i+1}-1}\sum_{f\in
{\bf Z},f^2>\left(t_1^2-4\right)\left(t_2^2-4\right)}h\left(t_1^2
-4,t_2^2-4,f\right)R_{t_1,t_2,f}\eqno (6.7)$$
with the abbreviations
$$a_{j_1,j_2}:=\left(-1\right)^{j_1+j_2}\left(\matrix{J\cr
j_1\cr}
\right)\left(\matrix{J\cr
j_2\cr}
\right),\eqno (6.8)$$
$$S_{t_1,t_2}:=\sum_{j_1=0}^J\sum_{j_2=0}^Ja_{j_1,j_2}{\cal J}\left(t_1,
t_2,k_{\left(X-j_1d\tau -2\right)/4},k_{\left(X-j_2d\tau -2\right
)/4}\right)\eqno (6.9)$$
and
$$R_{t_1,t_2,f}:=\sum_{j_1=0}^J\sum_{j_2=0}^Ja_{j_1,j_2}{\cal I}\left
(t_1,t_2,{f\over {\sqrt {t_1^2-4}\sqrt {t_2^2-4}}},k_{\left(X-j_1
d\tau -2\right)/4},k_{\left(X-j_2d\tau -2\right)/4}\right).\eqno
(6.10)$$
By (2.10) and Lemma 4.11 we see that the ${\cal I}$ and ${\cal J}$ functions involved in (6.9) and (6.10) can be nonzero only in the case
$$t_1^2-4\le X-j_1d\tau -2,\;t_2^2-4\le X-j_2d\tau -2.\eqno (6.11
)$$
If (6.11) is true, then by (2.10) and (4.1) we have that
$${\cal I}\left(t_1,t_2,{f\over {\sqrt {t_1^2-4}\sqrt {t_2^2-4}}}
,k_{\left(X-j_1d\tau -2\right)/4},k_{\left(X-j_2d\tau -2\right)/4}\right
)=Z\left(S_0,T_0,F\right),\eqno (6.12)$$
and by Lemma 4.11 we have that
$${\cal J}\left(t_1,t_2,k_{\left(X-j_1d\tau -2\right)/4},k_{\left(X-j_2d
\tau -2\right)/4}\right)=2\min\left(S_0,T_0\right),\eqno (6.13)$$
where we use the abbreviations
$$S_0=S_0\left(j_1,t_1\right)=\sqrt {{{X-j_1d\tau -2}\over {t_1^2
-4}}-1},\quad T_0=T_0\left(j_2,t_2\right)=\sqrt {{{X-j_2d\tau -2}\over {
t_2^2-4}}-1},\eqno (6.14)$$
$$F=F\left(t_1,t_2,f\right)=\left|{f\over {\sqrt {t_1^2-4}\sqrt {
t_2^2-4}}}\right|.\eqno (6.15)$$
We now consider the sum (6.7). Assume that (6.11) holds.
Then by Lemma 4.1 and by (4.8) we see that (6.12) can be nonzero only if

$$\left|f\right|<B\left(S_0,T_0\right)\sqrt {t_1^2-4}\sqrt {t_2^2
-4}.\eqno (6.16)$$
If (6.11) is true and (6.16) holds for some $f$ in (6.7), then
we have
$${{\left|f\right|}\over {\sqrt {t_1^2-4}\sqrt {t_2^2-4}}}-1\ge{1\over {\sqrt {
t_1^2-4}\sqrt {t_2^2-4}\left(\left|f\right|+\sqrt {t_1^2-4}\sqrt {
t_2^2-4}\right)}}\gg X^{-c}\eqno (6.17)$$
with some absolute constant $c$. If (6.11) holds, then we
determine (6.12) by Lemmas 4.1 and 4.3. In some cases we
will apply the upper bounds of Lemma 4.4. By (6.17) and
(4.12) we see that when we apply these lemmas for the estimation
of (6.12) we always have $\log{1\over {1-y}}\ll\log X$. So we get assuming (6.11) that for any $
f$ in (6.7) we have
that (6.12) is
$$\ll\left(\sqrt {{{X-j_1d\tau -2}\over {t_1^2-4}}-1}+\sqrt {{{X-
j_2d\tau -2}\over {t_2^2-4}}-1}\right)\log X.\eqno (6.18)$$
We note finally that if $a_i\le t_1,t_2<a_{i+1}$ for some $i$, then we have
$$\left(t_1t_2-5\right)^2+{{t_1t_2}\over 6}<\left(t_1^2-4\right)\left
(t_2^2-4\right)\le\left(t_1t_2-4\right)^2,\eqno (6.19)$$
since by the assumption $a_{i+1}\le{3\over 2}a_i$ made in (6.2) we have
${2\over 3}\le{{t_2}\over {t_1}}\le{3\over 2}$, and we also have $
t_1t_2\ge 9$. So there is an
absolute constant $c_0>0$ such that if $a_i\le t_1,t_2<a_{i+1}$, then
$$\left(t_1t_2-5\right)+c_0\le\sqrt {\left(t_1^2-4\right)\left(t_
2^2-4\right)}\le t_1t_2-4.\eqno (6.20)$$

{\bf 6.2. The case of very large} $a_i${\bf .} Assume that we have
$$\sqrt {X+2}-a_i=O\left({d\over {\sqrt X}}X^{\delta}\right)\eqno
(6.21)$$
for some $\delta >0$ which is chosen small enough in terms of
$\epsilon$. Consider first (6.7). Since it is easy to see that
$B\left(S,T\right)-1\le S^2+T^2$ for any $S,T\ge 0$, the number of
integers $f$ in (6.7) satisfying (6.16) is $\ll 1+\sqrt {t_1^2-4}\sqrt {
t_2^2-4}\left(\left({{X-j_1d\tau -2}\over {t_1^2-4}}-1\right)+\left
({{X-j_2d\tau -2}\over {t_2^2-4}}-1\right)\right)\ll dX^{\delta},$
where in the last step we used (6.21), (6.1) and (6.2). On the
other hand, we get by (6.11), (6.12), (6.18), (6.21) and (6.2) that
(6.10) is always $\ll_{\delta}{{\sqrt d}\over {\sqrt X}}X^{\delta}$ for every such $
f$. Hence for $i$ satisfying
(6.21) we have, applying also Lemma 3.1, (3.23) and (6.21) that (6.7) is
$$\ll_{\delta}X^{3\delta}{{d^{3/2}}\over {\sqrt X}}\sum_{t_1=a_i}^{
a_{i+1}-1}\sum_{t_2=a_i}^{a_{i+1}-1}S\left(t_1^2-4,t_2^2-4\right)
\ll_{\delta}X^{5\delta}d^{3/2}\left(\sqrt {{d\over {\sqrt X}}}+{{\left
({d\over {\sqrt X}}\right)^{{3\over 2}}}\over {X^{1/4}}}\right),$$
where we used (6.1), (6.2). By (6.3) we see that its contribution is acceptable in (5.18).

We now consider (6.6). We get by (6.11), (6.12), Lemma 4.1, (4.26), the
first relation in (4.27) and (6.21) that (6.10) is
$\ll_{\delta}{{X^{\delta}d}\over {\sqrt {\left(t_1^2-4\right)\left
(t_2^2-4\right)-f^2}}}.$ Hence for $i$ satisfying (6.21) we have, applying also Lemma
3.1 that (6.6) is
$$\ll_{\delta}X^{2\delta}d\sum_{t_1=a_i}^{a_{i+1}-1}\sum_{t_2=a_i}^{
a_{i+1}-1}\sum_{f\in {\bf Z},f^2<\left(t_1^2-4\right)\left(t_2^2-
4\right)}{{S\left(t_1^2-4,t_2^2-4,f^2\right)}\over {\sqrt {\left(
t_1^2-4\right)\left(t_2^2-4\right)-f^2}}}.\eqno (6.22)$$
In the $t_1\neq t_2$ part we use $S\left(t_1^2-4,t_2^2-4,f^2\right
)\le S\left(t_1^2-4,t_2^2-4\right)$, and we easily see that
the sum over $f$ in (6.22) is $\ll S\left(t_1^2-4,t_2^2-4\right)$. Then
applying (3.20) with $b=c$ and (6.21) we see that the $t_1\neq t_
2$
part of (6.22) is
$\ll_{\delta}X^{3\delta}d\sqrt {a_i}\left({d\over {\sqrt X}}X^{\delta}\right
)^{3/2}\ll_{\delta}X^{5\delta}{{d^{5/2}}\over {\sqrt X}}$, which is acceptable in
(5.18). For the $t_1=t_2$ part of (6.22) we use Lemma 3.7 and we get that the sum over $
f$
in (6.22) is $\ll_{\delta}X^{\delta}$, and so the  $t_1=t_2$ part of (6.22) is
$\ll_{\delta}X^{3\delta}d\left({d\over {\sqrt X}}X^{\delta}\right
)$, which is smaller than our estimate for the
$t_1\neq t_2$ part. So we proved that assuming (6.21) the
contribution of (6.6) is also acceptable in (5.18).

Assuming (6.21) in (6.5) we clearly have that
$\min\left(S_0,T_0\right)\ll_{\delta}X^{\delta}{{\sqrt d}\over {\sqrt
X}}$, hence using also (6.9), (6.13) and
Lemma 2.1 we get
that (6.5) is $\ll_{\delta}{{X^{2\delta}\sqrt d}\over {\sqrt X}}\left
({d\over {\sqrt X}}X^{\delta}\right)^2\sqrt X={{X^{4\delta}d^{5/2}}\over
X}.$ Hence in the case of (6.21) we get that the contribution
of $U_i$ in (5.18) is acceptable. We may assume from now
on that
$$\sqrt {X+2}-a_i>X^{\delta_0}{d\over {\sqrt X}}\eqno (6.23)$$
with some $\delta_0>0$, which is fixed in terms of $\epsilon$.

{\bf 6.3. Easy consequences of (6.23).} Observe that (6.23)
implies the relations
$${{X^{1/4}\sqrt {\sqrt X-a_i}}\over {a_i}}\ll S_0\left(j,t_1\right
),T_0\left(j,t_2\right)\ll{{X^{1/4}\sqrt {\sqrt X-a_i}}\over {a_i}}
,\eqno (6.24)$$
$${X\over {a_i^2}}\ll 1+S_0^2\left(j,t_1\right),1+T_0^2\left(j,t_
2\right)\ll{X\over {a_i^2}},\eqno (6.25)$$
$$T_0\left(j_1,t\right)-T_0\left(j_2,t\right)=O\left({d\over {X^{
1/4}a_i\sqrt {\sqrt X-a_i}}}\right)=S_0\left(j_1,t\right)-S_0\left
(j_2,t\right)\eqno (6.26)$$
and
$$S_0\left(j_1,t_1\right)-T_0\left(j_2,t_2\right)={{X\left(t_2^2-
t_1^2\right)\left(\left(t_1^2-4\right)\left(t_2^2-4\right)\right)^{
-1}}\over {S_0\left(j_1,t_1\right)+T_0\left(j_2,t_2\right)}}+O\left
({{dX^{-1/4}}\over {\sqrt {\sqrt X-a_i}a_i}}\right)\eqno (6.27)$$
for every real numbers $0\le j_1,j_2,j\le J$ and integers $a_i\le
t_1,t_2,t\le a_{i+1}$.

We have in general that
${d\over {dT}}\left(\sqrt {1+S^2}\sqrt {1+T^2}+ST\right)={{\sqrt {
1+S^2}T}\over {\sqrt {1+T^2}}}+S$ and
$${d\over {dT}}\left(\sqrt {1+S^2}\sqrt {1+T^2}-ST\right)={{\sqrt {
1+S^2}T}\over {\sqrt {1+T^2}}}-S={{\left(T-S\right)\left(T+S\right
)}\over {\sqrt {1+T^2}\left(\sqrt {1+S^2}T+\sqrt {1+T^2}S\right)}}
,$$
hence we get by (6.24)-(6.27) and the mean-value theorem that
$$B\left(S_0\left(j_1,t_1\right),T_0\left(j,t_2\right)\right)-B\left
(S_0\left(j_1,t_1\right),T_0\left(0,t_2\right)\right)\ll{d\over {
a_i^2}}\eqno (6.28)$$
and
$$A\left(S_0\left(j_1,t_1\right),T_0\left(j,t_2\right)\right)-A\left
(S_0\left(j_1,t_1\right),T_0\left(0,t_2\right)\right)\ll\left({{\left
|t_2-t_1\right|}\over {a_i}}+{d\over X}\right){d\over {X-a_i^2}}\eqno
(6.29)$$
for every real numbers $0\le j_1,j\le J$ and integers $a_i\le t_1
,t_2\le a_{i+1}$.

We will also need later the easily proved general identity
$${{T_0}\over {S_0}}-A\left(S_0,T_0\right)={{\left(1+S_0^2\right)\left
(T_0^2-S_0^2\right)}\over {S_0T_0\left(1+S_0^2\right)+S_0^2\sqrt {\left
(1+S_0^2\right)\left(1+T_0^2\right)}}}.\eqno (6.30)$$
This easily implies by (6.24) and (6.25) that assuming (6.23) we
have
$$\left|{{T_0}\over {S_0}}-1\right|\ll\left|{{T_0}\over {S_0}}-A\left
(S_0,T_0\right)\right|\ll\left|{{T_0}\over {S_0}}-1\right|\eqno (
6.31)$$
for every choice $S_0=S_0\left(j_1,t_1\right)$, $T_0=T_0\left(j_2
,t_2\right)$ with any
real numbers $0\le j_1,j_2\le J$ and integers $a_i\le t_1,t_2\le
a_{i+1}$.
We also see from (6.30) that
the signs of ${{T_0}\over {S_0}}-A\left(S_0,T_0\right)$ and ${{T_
0}\over {S_0}}-1$ are the same.
Therefore we get from (6.31) that assuming (6.23) we
have
$$\left|{{T_0\left(j_2,t_2\right)}\over {S_0\left(j_1,t_1\right)}}
-1\right|\ll\left|{{T_0\left(j_2,t_2\right)}\over {S_0\left(j_1,t_
1\right)}}-F\right|\ll\left|{{T_0\left(j_2,t_2\right)}\over {S_0\left
(j_1,t_1\right)}}-1\right|\eqno (6.32)$$
for any real numbers $0\le j_1,j_2\le J$ and for any
$1<F\le A\left(S_0\left(j_1,t_1\right),T_0\left(j_2,t_2\right)\right
)$.

Assuming (6.23) we see by (6.2) that (6.11) is always true.
It follows then by (6.9) and (6.13) that (6.5) equals
$$2\sum_{t_1=a_i}^{a_{i+1}-1}\sum_{t_2=a_i}^{a_{i+1}-1}E_{t_1,t_2}
\sum_{j_1=0}^J\sum_{j_2=0}^Ja_{j_1,j_2}\min\left(S_0\left(j_1,t_1\right
),T_0\left(j_2,t_2\right)\right).\eqno (6.33)$$
We also see that (6.12) always holds. Then by Lemma 4.1
and (6.10) we get for any $f$ that
$$R_{t_1,t_2,f}=2\sum_{j_1=0}^J\sum_{j_2=0}^Ja_{j_1,j_2}\left(J\left
(S_0,T_0,F\right)+J\left(T_0,S_0,F\right)\right),\eqno (6.34)$$
By the change of variables $j_1\rightarrow j_2$, $t_1\rightarrow
t_2$ we see that
substituting (6.34) into (6.7) and (6.6) the contributions of
$J\left(S_0,T_0,F\right)$ and $J\left(T_0,S_0,F\right)$ in (6.7) are the same.

{\bf 6.4. Estimating (6.5).} Assume besides (6.23) that we have
$$\left|t_2-t_1\right|>{{da_i}\over X}X^{\delta}\eqno (6.35)$$
for some $\delta >0$ which is fixed in terms of $\epsilon$. Then we
see from (6.27) and (6.24) that the sign of
$S_0\left(j_1,t_1\right)-T_0\left(j_2,t_2\right)$ is the same for every pair
$0\le j_1,j_2\le J$. But then that part of (6.33) where (6.35)
holds is $0$, since $\sum_{j_1=0}^Ja_{j_1,j_2}=0$ for every $j_2$, and
$\sum_{j_2=0}^Ja_{j_1,j_2}=0$ for every $j_1$ by (6.8). So we may assume
in (6.33) that
$$\left|t_2-t_1\right|\ll{{da_i}\over X}X^{\delta}\eqno (6.36)$$
for some $\delta >0$ which is chosen small enough in terms of
$\epsilon$. Then we see using (6.24) and Lemma 2.1 that (6.5) is $
\ll_{\delta}X^{2\delta}X^{1/4}\sqrt {\sqrt X-a_i}a_i\left(1+{{da_
i}\over X}\right)\ll_{\delta}X^{2\delta}\sqrt Xd,$
which is acceptable in (5.18).

{\bf 6.5. Estimating (6.6).} Assume besides (6.23) that in (6.6) we have
(6.36) and
$$1-\left({d\over {X-a_i^2}}\right)^2X^{\delta}<F<1\eqno (6.37)$$
for some $\delta >0$ which is chosen small enough in terms of
$\epsilon$. By (6.34) and (4.26) we have
$$R_{t_1,t_2,f}=2\sum_{j_1=0}^J\sum_{j_2=0}^Ja_{j_1,j_2}\left(K\left
(S_0,T_0,F\right)+K\left(T_0,S_0,F\right)\right)\eqno (6.38)$$
in the case $F<1$, and by the substitutions $j_1\rightarrow j_2$, $
t_1\rightarrow t_2$ we see that the
contributions of $K\left(S_0,T_0,F\right)$ and $K\left(T_0,S_0,F\right
)$ in (6.6) are
the same. Hence it is enough to consider the contribution
of $K\left(S_0,T_0,F\right)$. Applying the second relation in (4.27) we see for fixed $
t_1$, $t_2$, $f$ and $j_1$ that
$$\sum_{j_2=0}^J\left(-1\right)^{j_2}\left(\matrix{J\cr
j_2\cr}
\right)K\left(S_0,T_0\left(j_2,t_2\right),F\right)\ll\max_{0\le j_
2<J}{{\left|S_0\right|\left|T_0\left(j_2+1,t_2\right)-T_0\left(j_
2,t_2\right)\right|}\over {\sqrt {1-F^2}}}.$$
The parameters are written here only in the case of $T_0$,
since only this variable depends on $j_2$. By (6.24) and
(6.26) this is $\ll{d\over {a_i^2\sqrt {1-F^2}}}$. Hence using (6.8), (6.38) and Lemma 3.1 we
get that that part of (6.6) where (6.36) and (6.37) hold is
$$\ll_{\delta}{{dX^{\delta}}\over {a_i^2}}\sum_{t_1=a_i}^{a_{i+1}
-1}\sum_{a_i\le t_2<a_{i+1},\left|t_2-t_1\right|\ll da_iX^{\delta
-1}}\sum_{f\in {\bf Z},0<1-F<\left({d\over {X-a_i^2}}\right)^2X^{
\delta}}{{S\left(t_1^2-4,t_2^2-4,f^2\right)}\over {\sqrt {1-F^2}}}
.\eqno (6.39)$$
By (6.20) we see that the sum over $f$ here is
$$\ll a_i\sum_{t_1t_2-5-\sqrt {\left(t_1^2-4\right)\left(t_2^2-4\right
)}\left({d\over {X-a_i^2}}\right)^2X^{\delta}\le\left|f\right|\le
t_1t_2-5}{{S\left(t_1^2-4,t_2^2-4,f^2\right)}\over {\sqrt {\left(
t_1t_2-4\right)-\left|f\right|}}}.\eqno (6.40)$$
In the case $t_1\neq t_2$ we use
$S\left(t_1^2-4,t_2^2-4,f^2\right)\le S\left(t_1^2-4,t_2^2-4\right
)$, and we see that (6.40)
is $\ll_{\delta}X^{\delta}a_i\left(1+{{da_i}\over {X-a_i^2}}\right
)S\left(t_1^2-4,t_2^2-4\right).$ So, applying also
(3.20) we get that the
$t_1\neq t_2$ part of (6.39) is
$$\ll_{\delta}X^{3\delta}\left({d\over {a_i}}+{{d^2}\over {X-a_i^
2}}\right)\left(\left(\min\left(a_i,\sqrt X-a_i\right)\right)\left
(a_i{{da_i}\over X}\right)^{1/2}\right)\ll_{\delta}X^{3\delta}{{d^{
5/2}}\over {\sqrt X}},$$
which is acceptable in (5.18). In the case $t_1=t_2$ we
estimate (6.40) by Lemma 3.7 and we get that (6.40) is
$\ll_{\delta}X^{\delta}a_i\left(1+\left(t_1^2-4\right)\left({d\over {
X-a_i^2}}\right)^2\right)^{1/2}.$ So the $t_1=t_2$ part of (6.39) is
$$\ll_{\delta}X^{2\delta}{d\over {a_i}}\left(\min\left(a_i,\sqrt
X-a_i\right)\right)\left(1+a_i{d\over {X-a_i^2}}\right)\ll_{\delta}
X^{2\delta}{{d^2}\over {\sqrt X}}.$$
Hence that part of (6.6) where (6.36) and (6.37) hold is acceptable in (5.18).

So it is enough to consider that part of (6.6) where at least one of the
conditions (6.36) and (6.37) is false. We prove that this part is negligible. We
use (6.34), and we recall that the contributions of
$J\left(S_0,T_0,F\right)$ and $J\left(T_0,S_0,F\right)$ in (6.6) are the same. By
Lemma 4.8 and $\sum_{j_2=0}^Ja_{j_1,j_2}=0$ we see that it is enough
to show that for fixed $t_1$, $t_2$, $f$, $j_1$ the sum
$$\sum_{j_2=0}^J\left(-1\right)^{j_2}\left(\matrix{J\cr
j_2\cr}
\right)V\left(s_i\left(S_0\left(j_1,t_1\right),T_0\left(j_2,t_2\right
),F\left(t_1,t_2,f\right)\right)\right)\eqno (6.41)$$
is negligibly small for $i=1,2$. Observe that by the
notation of Lemma 4.10, using $t=t_2$, $S_0=S_0\left(j_1,t_1\right
)$, $F=F\left(t_1,t_2,f\right)$ and $\tau =1$ for
$i=1$, $\tau =-1$ for $i=2$ we have that
$$V\left(s_i\left(S_0\left(j_1,t_1\right),T_0\left(j_2,t_2\right)
,F\left(t_1,t_2,f\right)\right)\right)=k\left(X-j_2d\tau -2\right
).$$
Then by Taylor's formula with remainder (see e.g.
Theorem 7.6 of [A]) we have that (6.41) is
$\ll d^J\max_{X-2-Jd\tau\le x\le X-2}\left|k^{\left(J\right)}\left
(x\right)\right|.$ Then by Lemma 4.10 and
(6.24) we see that (6.41) is
$$\ll\sqrt {1+{1\over {S_0\left(j_1,t_1\right)^2}}}d^J\left(X-a^2_
i\right)^{-J}\max\left(1,\left({1\over {\left(\sqrt {1-F^2}+\left
|\tau F+{{T_0\left(j,t_2\right)}\over {S_0\left(j_1,t_1\right)}}\right
|\right)}}\right)^J\right)\eqno (6.42)$$
with some real number $0\le j\le J$. We see that if (6.37) is false, then this is
negligibly small, since $J$ is fixed to be large enough. So
we can assume that (6.37) is true but (6.36) is false. We show
that (6.42) is negligibly small. If
$$\left|\tau F+{{T_0\left(j,t_2\right)}\over {S_0\left(j_1,t_1\right
)}}\right|\gg X^{\delta}{d\over {X-a_i^2}},\eqno (6.43)$$
then this is true. So we may assume that (6.43) is false.
But then using also (6.37) and the triangle inequality,
taking into account (6.23) we get
$$\left|\tau +{{T_0\left(j,t_2\right)}\over {S_0\left(j_1,t_1\right
)}}\right|\ll X^{\delta}{d\over {X-a_i^2}}.\eqno (6.44)$$
This is impossible for $\tau =1$ for small $\delta$ by (6.23), so we may assume $
\tau =-1$. Hence,
using (6.27) and (6.44) with $\tau =-1$, applying also (6.24) we
get $\left|t_2-t_1\right|\ll X^{\delta}{}^{}{{da_i}\over X}.$ But this is a contradiction, since we assumed that (6.36)
is false. So that part of (6.6) where at least one of the
conditions (6.36) and
(6.37) is false is also negligibly small. Consequently (6.6) is acceptable in (5.18).

{\bf 6.6. A new expression for (}6.7{\bf ).} Recall that substituting
(6.34) into (6.7) the contributions of $J\left(S_0,T_0,F\right)$ and
$J\left(T_0,S_0,F\right)$ in (6.7)
are the same. Hence, applying also Lemma 4.3 we get that
(6.7) equals
$$4\sum_{t_1=a_i}^{a_{i+1}-1}\sum_{t_2=a_i}^{a_{i+1}-1}\left(A_{t_
1,t_2}+B_{t_1,t_2}-C_{t_1,t_2}\right)\eqno (6.45)$$
with
$$A_{t_1,t_2}:=\sum_{0\le j_1,j_2\le J,f\in {\bf Z},1<F\le B\left
(S_0,T_0\right)}h\left(t_1^2-4,t_2^2-4,f\right)a_{j_1,j_2}S_0\Phi\left
(y_1\left(S_0,T_0,F\right)\right),\eqno (6.46)$$
$$B_{t_1,t_2}:=\sum_{0\le j_1,j_2\le J,f\in {\bf Z},T_0\ge S_0,1<
F\le A\left(S_0,T_0\right)}h\left(t_1^2-4,t_2^2-4,f\right)a_{j_1,
j_2}S_0\Phi\left(y_2\left(S_0,T_0,F\right)\right),\eqno (6.47)$$
$$C_{t_1,t_2}:=\sum_{0\le j_1,j_2\le J,f\in {\bf Z},T_0\le S_0,1<
F\le A\left(S_0,T_0\right)}h\left(t_1^2-4,t_2^2-4,f\right)a_{j_1,
j_2}S_0\Phi\left(y_2\left(S_0,T_0,F\right)\right).\eqno (6.48)$$
{\bf 6.7. The contribution of} $B_{t_1,t_2}$ {\bf and} $C_{t_1,t_
2}${\bf .} Assume besides
(6.23) that
$$1<\left|F\right|\le A\left(S_0,T_0\right)\eqno (6.49)$$
and (6.36) holds with some $\delta >0$ which is chosen small enough in terms of
$\epsilon$. Applying (6.17), (4.12) and (4.16) we get that the terms
$\Phi\left(y_2\left(S_0,T_0,F\right)\right)$ in (6.47), (6.48) are always $
O\left(\log X\right)$. Using
$$A\left(S_0,T_0\right)-1={{\left(S_0-T_0\right)^2}\over {\sqrt {\left
(1+S_0^2\right)\left(1+T_0^2\right)}+S_0T_0+1}}\le{{\left(S_0-T_0\right
)^2}\over {\sqrt {\left(1+S_0^2\right)\left(1+T_0^2\right)}}},$$
(6.27), (6.25), (6.24) and (6.36) we see that the number of integers $
f$
satisfying (6.49) and (6.15) is $\ll 1+X^{2\delta}{{d^2a_i^2/X}\over {
X-a^2_i}}.$ So we get, applying also Lemma
3.1 that the contribution to (6.45) of the terms $B_{t_1,t_2}$,
$C_{t_1,t_2}$ satisfying (6.49) and (6.36) is
$$\ll_{\delta}X^{3\delta}\left({{\sqrt {X-a^2_i}}\over {a_i}}+{{d^
2a_i/X}\over {\sqrt {X-a^2_i}}}\right)\sum_{t_1=a_i}^{a_{i+1}-1}\sum_{
a_i\le t_2\le a_{i+1}-1,\left|t_2-t_1\right|\ll a_i{d\over X}X^{\delta}}
S\left(t_1^2-4,t_2^2-4\right).$$
By (3.20) and (3.22) it is
$$\ll_{\delta}X^{4\delta}\left({{\sqrt {X-a^2_i}}\over {a_i}}+{{d^
2a_i/X}\over {\sqrt {X-a^2_i}}}\right)\left(a_i\sqrt {\sqrt X-a_i}
+\left(\sqrt X-a_i\right)a_i\left({d\over X}\right)^{1/2}\right),$$
and this is $\ll_{\delta}X^{4\delta}{{d^{5/2}}\over {\sqrt X}}$ by (6.1), (6.2) and (5.19). Hence we
proved that the contribution to (6.45) of the terms $B_{t_1,t_2}$, $
C_{t_1,t_2}$ satisfying (6.49) and (6.36) is acceptable in (5.18).

Consider now the contribution of those terms
$B_{t_1,t_2}$$,C_{t_1,t_2}$ to (6.45) for which
$$1<F\le A\left(S_0\left(j_1,t_1\right),T_0\left(0,t_2\right)\right
)-{{d\left|t_2-t_1\right|}\over {a_i\left(X-a_i^2\right)}}X^{\delta}\eqno
(6.50)$$
and (6.35) hold for some $\delta >0$ which is fixed in terms of
$\epsilon$. We want to prove that this contribution is negligibly
small. We first show that for fixed $t_1,t_2,j_1$ and $f$ the conditions in
the summations in (6.47) and (6.48) are independent of
$0\le j_2\le J$. It is enough to see that we have
$F\le A\left(S_0\left(j_1,t_1\right),T_0\left(j,t_2\right)\right)$ for every $
0\le j\le J$, and the sign of $S_0\left(j_1,t_1\right)-T_0\left(j
,t_2\right)$
is the same for every $0\le j\le J$. These statements follow
easily from (6.29), (6.27) and (6.24). Hence for fixed
$t_1,t_2,j_1$ and $f$ satisfying $a_i\le t_1,t_2<a_{i+1}$, $0\le
j_1\le J$ and the
conditions (6.50), (6.35) we have that either each $0\le j_2\le J$ satisfies the
conditions of the summations in (6.47), or
each $0\le j_2\le J$ satisfies the conditions of the summations
in (6.48). Consequently, recalling (6.8) we see that it is
enough to show that for every fixed $t_1,t_2,j_1$ and $f$ satisfying the above-mentioned
conditions the sum
$$\sum_{j_2=0}^J\left(-1\right)^{j_2}\left(\matrix{J\cr
j_2\cr}
\right)\Phi\left(y_2\left(S_0\left(j_1,t_1\right),T_0\left(j_2,t_
2\right),F\left(t_1,t_2,f\right)\right)\right)\eqno (6.51)$$
is negligibly small. Observe that by the notation of Lemma
4.6, using $\tau =-1$ and $t=t_2$, $S_0=S_0\left(j_1,t_1\right)$, $
F=F\left(t_1,t_2,f\right)$
there we have
$$\Phi\left(y_2\left(S_0\left(j_1,t_1\right),T_0\left(j_2,t_2\right
),F\left(t_1,t_2,f\right)\right)\right)=K\left(X-j_2d\tau -2\right
).$$
Theorem 7.6 of [A] shows that (6.51) is $\ll d^J\max_{X--2-Jd\tau
\le x\le X-2}\left|K^{\left(J\right)}\left(x\right)\right|.$ By Lemma 4.6, (6.32),
(6.27), (6.24), (6.23), (6.1), (6.2) this is
$$\ll d^J\left(X-a^2_i\right)^{-J}\max\left(1,\left({{\left|t_2-t_
1\right|}\over {a_i\left(A\left(S_0,T_0\right)-F\right)}}\right)^
J\right).$$
By (6.50), (6.29) and (6.23) we see that this is negligibly
small, since $J$ is fixed to be large enough in terms of
$\epsilon$. Hence we proved that the contribution to (6.45) of those terms
$B_{t_1,t_2}$$,C_{t_1,t_2}$ for which (6.50) and (6.35) hold is negligibly small.

Consider the contribution to (6.45) of those terms $B_{t_1,t_2}$,
$C_{t_1,t_2}$ for which (6.35) holds and we have
$$A\left(S_0\left(j_1,t_1\right),T_0\left(0,t_2\right)\right)-{{d\left
|t_2-t_1\right|}\over {a_i\left(X-a_i^2\right)}}X^{\delta}<F\le A\left
(S_0\left(j_1,t_1\right),T_0\left(j_2,t_2\right)\right)\eqno (6.5
2)$$
for some $\delta >0$ which is small enough in terms of
$\epsilon$. Using also (6.29) this shows that the number of possible
values of the integers $f$ satisfying (6.52) is
$$\ll 1+{{d\left|t_2-t_1\right|a_i}\over {\left(X-a_i^2\right)}}X^{
\delta}.\eqno (6.53)$$
It is easy to compute that
$$y_2\left(S_0,T_0,F\right)=\sqrt {{{\left(B\left(S_0,T_0\right)+
F\right)\left(A\left(S_0,T_0\right)-F\right)}\over {\left(T_0-S_0
F\right)^2}}},$$
so using (6.24), (6.25), (6.27), (6.32) and (6.52) we get
$$y_2\left(S_0,T_0,F\right)\ll_{\delta}X^{\delta}\sqrt {{{da_i}\over {
X\left|t_2-t_1\right|}}},\;S_0y_2\left(S_0,T_0,F\right)\ll_{\delta}
X^{\delta}\sqrt {{{d\left(\sqrt X-a_i\right)}\over {a_i\sqrt X\left
|t_2-t_1\right|}}}.\eqno (6.54)$$
Now, if $a_i\ge{{\sqrt X}\over 2}$, then we have $S_0\ll 1$ by (6.24), and so
by Lemma 4.4, (4.12) and (6.17) we get that
$$S_0\Phi\left(y_2\left(S_0,T_0,F\right)\right)\ll_{\delta}S_0y_2^
3\left(S_0,T_0,F\right)X^{\delta}\ll_{\delta}X^{4\delta}\sqrt {{{
d^3\left(\sqrt X-a_i\right)}\over {X^2\left|t_2-t_1\right|^3}}}.\eqno
(6.55)$$
If $a_i\le{{\sqrt X}\over 2}$, then we have  $S_0\gg 1$ by (6.24),
and by the second relation in (6.54), Lemma 4.4, (4.12), (6.17) we
see that if $\left|t_2-t_1\right|\ll{d\over {a_i}}$, then
$$S_0\Phi\left(y_2\left(S_0,T_0,F\right)\right)\ll_{\delta}X^{\delta}
S_0y_2\left(S_0,T_0,F\right)\ll_{\delta}X^{2\delta}\sqrt {{d\over {
a_i\left|t_2-t_1\right|}}},\eqno (6.56)$$
while if $\left|t_2-t_1\right|\gg{d\over {a_i}}$, then
$$S_0\Phi\left(y_2\left(S_0,T_0,F\right)\right)\ll_{\delta}X^{\delta}
S_0^3y_2^3\left(S_0,T_0,F\right)\ll_{\delta}X^{4\delta}\sqrt {{{d^
3}\over {a_i^3\left|t_2-t_1\right|^3}}}.\eqno (6.57)$$
If $a_i\ge{{\sqrt X}\over 2}$, then by (6.35) and (5.19) we see that the
second term is larger than the first one in (6.53). Then by
(6.53) and (6.55) we see that the contribution to (6.45) of
those terms $B_{t_1,t_2}$, $C_{t_1,t_2}$ for which (6.35) and (6.52) hold is
$$\ll_{\delta}X^{6\delta}{{d^{5/2}}\over {X^{3/4}\sqrt {X-a_i^2}}}
\sum_{t_1=a_i}^{a_{i+1}-1}\sum_{a_i\le t_2<a_{i+1},\left|t_2-t_1\right
|\ge da_iX^{\delta -1}}{{S\left(t^2_1-4,t^2_2-4\right)}\over {\sqrt {\left
|t_2-t_1\right|}}}\eqno (6.58)$$
in the case $a_i\ge{{\sqrt X}\over 2}$.  By (3.21) we have that the sum
over $t_1$,$t_2$ is $\ll_{\delta}X^{\delta}\left(\sqrt X-a_i\right
)X^{1/4}$, hence (6.58) is acceptable in (5.18).

If $a_i\le{{\sqrt X}\over 2}$, then by (6.53) and (6.56) we see that the contribution to (6.45) of
those terms $B_{t_1,t_2}$, $C_{t_1,t_2}$ for which (6.35), (6.52) and $\left
|t_2-t_1\right|\ll{d\over {a_i}}$ hold is $\ll_{\delta}$ than the
sum of
$$X^{4\delta}{{d^2}\over X}\sum_{t_1=a_i}^{a_{i+1}-1}\sum_{a_i\le
t_2<a_{i+1},0<\left|t_2-t_1\right|\ll{d\over {a_i}}}S\left(t^2_1-
4,t^2_2-4\right)\eqno (6.59)$$
and
$$X^{4\delta}\sqrt {{d\over {a_i}}}\sum_{t_1=a_i}^{a_{i+1}-1}\sum_{
a_i\le t_2<a_{i+1},0<\left|t_2-t_1\right|}{{S\left(t^2_1-4,t^2_2-
4\right)}\over {\sqrt {\left|t_2-t_1\right|}}}.\eqno (6.60)$$
By (3.20) we get that (6.59) is
$\ll_{\delta}X^{5\delta -1}d^2a_i^{3/2}\left(d/a_i\right)^{1/2}\ll_{
\delta}X^{5\delta}d^{5/2}/\sqrt X$$.$ We estimate (6.60) by
(3.21), the result is an upper bound $X^{5\delta}\left(d/a_i\right
)^{1/2}a_i^{3/2}\ll X^{5\delta}\sqrt d\sqrt X$, which is smaller than $
d^{5/2}/\sqrt X$
by (5.19).

If $a_i\le{{\sqrt X}\over 2}$ and $\left|t_2-t_1\right|\gg{d\over {
a_i}}$, then by (5.19) we see that
the second term in (6.53) is larger than the first one.
Hence by (6.53) and (6.57) we see in the case $a_i\le{{\sqrt X}\over
2}$ that the contribution to (6.45) of
those terms $B_{t_1,t_2}$, $C_{t_1,t_2}$ for which (6.35), (6.52) and
$\left|t_2-t_1\right|\gg{d\over {a_i}}$ hold is
$$\ll_{\delta}X^{6\delta}{{d^{5/2}}\over {a_i^{1/2}X}}\sum_{t_1=a_
i}^{a_{i+1}-1}\sum_{a_i\le t_2<a_{i+1},0<\left|t_2-t_1\right|}{{S\left
(t^2_1-4,t^2_2-4\right)}\over {\sqrt {\left|t_2-t_1\right|}}}\ll_{
\delta}X^{7\delta}{{d^{5/2}a_i}\over X},$$
where in the last step we used (3.21). This is again
acceptable in (5.18).

We examined every case, so we proved that
the contribution to (6.45) of those terms $B_{t_1,t_2}$, $C_{t_1,
t_2}$ for which (6.35) and (6.52) hold is acceptable
in (5.18). Using also the previous estimates we see that the whole contributions
of $B_{t_1,t_2}$ and $C_{t_1,t_2}$ in (6.45) is acceptable in (5.18).

{\bf 6.8. The contribution of} $A_{t_1,t_2}${\bf .} Consider now that part of
the contribution of $A_{t_1,t_2}$  in (6.45) where
$$1<F\le B\left(S_0\left(j_1,t_1\right),T_0\left(0,t_2\right)\right
)-{d\over {a_i^2}}X^{\delta}\eqno (6.61)$$
for some $\delta >0$ which is fixed in terms of $\epsilon$. We show that for fixed $
t_1,t_2,j_1$ and $f$ the condition
in the summation in (6.46) is independent of $0\le j_2\le J$. It
is enough to see that we have $F\le B\left(S_0\left(j_1,t_1\right
),T_0\left(j,t_2\right)\right)$
for every $0\le j\le J$, and this follows by (6.28). Hence for fixed $
t_1,t_2,j_1$ and $f$ satisfying
$a_i\le t_1,t_2<a_{i+1}$, $0\le j_1\le J$ and (6.61) each $0\le j_
2\le J$ satisfies the
conditions of the summation in (6.46). Consequently,
recalling (6.8) we see that if we can show that for
every fixed $t_1,t_2,j_1$ and $f$ satisfying the above-mentioned
conditions the sum
$$\sum_{j_2=0}^J\left(-1\right)^{j_2}\left(\matrix{J\cr
j_2\cr}
\right)\Phi\left(y_1\left(S_0\left(j_1,t_1\right),T_0\left(j_2,t_
2\right),F\left(t_1,t_2,f\right)\right)\right)\eqno (6.62)$$
is negligibly small, then we will get that that part of the
contribution of $A_{t_1,t_2}$ in (6.45) where (6.61) is true is
negligibly small. Observe that by the notations of Lemma 4.6, using $
\tau =1$ and $t=t_2$, $S_0=S_0\left(j_1,t_1\right)$, $F=F\left(t_
1,t_2,f\right)$
there we have
$$\Phi\left(y_1\left(S_0\left(j_1,t_1\right),T_0\left(j_2,t_2\right
),F\left(t_1,t_2,f\right)\right)\right)=K\left(X-j_2d\tau -2\right
).$$
Theorem 7.6 of [A] gives that (6.62) is
$\ll d^J\max_{X-Jd\tau -2\le x\le X-2}\left|K^{\left(J\right)}\left
(x\right)\right|.$ By Lemma 4.6 and
(6.24) this is $\ll d^J\max\left(\left(X-a^2_i\right)^{-J},\left({
1\over {a_i^2(\left(B\left(S_0,T_0\right)-F\right)}}\right)^J\right
).$
By (6.61), (6.28) and (6.23) this is negligibly small, since $J$
is fixed to be large enough in terms of $\epsilon$. Hence that part of the contribution of $
A_{t_1,t_2}$ in (6.45) where (6.61) holds is negligibly small.

Consider now that part of the contribution of $A_{t_1,t_2}$ in (6.45) where
$$B\left(S_0\left(j_1,t_1\right),T_0\left(0,t_2\right)\right)-{d\over {
a_i^2}}X^{\delta}<F\le B\left(S_0\left(j_1,t_1\right),T_0\left(j_
2,t_2\right)\right)\eqno (6.63)$$
for some $\delta >0$ which is chosen small enough in terms of
$\epsilon$. It is easy to compute that
$$y_1\left(S_0,T_0,F\right)=\sqrt {{{\left(B\left(S_0,T_0\right)-
F\right)\left(A\left(S_0,T_0\right)+F\right)}\over {\left(T_0+S_0
F\right)^2}}}.$$
By (6.25), (6.63) and (5.19) we see ${X\over {a_i^2}}\ll B\left(S_
0,T_0\right)\ll{X\over {a_i^2}}$ and $B\left(S_0,T_0\right)-F=o\left
({X\over {a_i^2}}\right)$. So
$$S_0y_1\left(S_0,T_0,F\right)\ll{{\sqrt {B\left(S_0,T_0\right)-F}}\over {\sqrt {
B\left(S_0,T_0\right)}}}=o\left(1\right).\eqno (6.64)$$
Then applying Lemma 4.4 (note that (4.14) cannot hold by (6.64)), using also
(4.12) and (6.17) we get in every case that
$$S_0\Phi\left(S_0,y_1\left(S_0,T_0,F\right)\right)\ll_{\delta}X^{
\delta}\left(S_0+S_0^3\right)y_1^3\left(S_0,T_0,F\right)\ll_{\delta}{{
X^{1+\delta}\left(B\left(S_0,T_0\right)-F\right)^{3/2}}\over {a_i^
2S_0^2\left(B\left(S_0,T_0\right)\right)^{3/2}}},$$
the second inequality follows by (6.25) and (6.64). By (6.24), (6.25) and
(6.63) this gives
$$S_0\Phi\left(S_0,y_1\left(S_0,T_0,F\right)\right)\ll_{\delta}X^{
3\delta}{{d^{3/2}}\over {X\left(\sqrt X-a_i\right)}}.$$
The number of possible values of the integers $f$
satisfying (6.63) is $\ll_{\delta}X^{\delta}d$, therefore using also Lemma 3.1
we get that that part of the
contribution of $A_{t_1,t_2}$ in (6.45) where (6.63) holds is
$$\ll_{\delta}X^{5\delta}{{d^{5/2}}\over {X\left(\sqrt X-a_i\right
)}}\sum_{t_1=a_i}^{a_{i+1}-1}\sum_{t_2=a_i}^{a_{i+1}-1}S\left(t^2_
1-4,t^2_2-4\right)\ll_{\delta}{{X^{6\delta}d^{5/2}}\over {X\left(\sqrt
X-a_i\right)}}\left(\sqrt X-a_i\right)a_i,$$
where in the last step we used (3.23), noting that $a\left(b-a\right
)$
is an upper bound there for both terms. This estimate is
again acceptable in (5.18). So we proved that the
contribution of $A_{t_1,t_2}$ in (6.45) is acceptable in (5.18),
hence the whole sum (6.45) is acceptable. The proof of
(5.18) is now complete, so Theorem 1.1 is also proved.

 \bigskip\noindent {\bf References}

\nobreak
\parindent=12pt
\nobreak

\item{[A]} T. Apostol, {Calculus, Vol. 1, Second Edition},
{\it Wiley}, 1967

\item{[B1]} A. Bir\'o, {\it Local average of the hyperbolic circle problem for Fuchsian groups}, Mathematika 64 (2018), 159--183.

\item{[B2]} A. Bir\'o, {\it Cycle integrals of Maass forms of weight 0 and Fourier coefficients of Maass forms of weight 1/2}, Acta Arithmetica,
94 (2) (2000), 103-152.

\item{[Bu]} D. A. Buell, {Binary quadratic forms: classical theory and modern computations,}

 {\it Springer}, 1989

\item{[C]} F. Chamizo, {\it Some applications of large sieve in Riemann surfaces}, Acta Arithmetica, 77 (4) (1996), 315-337.

\item{[C-R]} G. Cherubini, M. Risager, {\it On the variance of the error term in the hyperbolic circle problem}, Revista matem\'atica iberoamericana, 34 (2) (2018), 655-685.

\item{[H-W]} K. Hardy, S. Williams, {\it The class number of pairs of positive-definite binary quadratic forms}, Acta Arithmetica, 52 (2) (1989), 103-117.

\item{[I]} H. Iwaniec, {Introduction to the spectral theory of automorphic forms,} {\it Rev. Mat. Iberoamericana}, 1995

\item{[L-P]} P. Lax, R. Phillips {\it The asymptotic distribution of lattice points in Euclidean and non-Euclidean spaces}, J. Funct. Anal. 46 (1982), 280-350.

\item{[M]} J. Morales, {\it The classification of pairs of binary quadratic forms}, Acta Arithmetica. 59 (2) (1991), 105-121.

\item{[P-R]} Y. Petridis, M. Risager, {\it Local average in hyperbolic lattice point counting, with an appendix by Niko Laaksonen}, Math. Zeitschrift, 285 (2017), 1319-1344.

\bye